\documentclass[12pt]{amsart}
\usepackage{graphicx}
\usepackage{amsfonts}
\usepackage{amsmath}
\usepackage{amsthm}
\usepackage{amssymb}
\usepackage{geometry}
\usepackage{thmtools}
\usepackage{tikz}
\usetikzlibrary{calc,trees,positioning,arrows,fit,shapes,calc}
\usepackage{circuitikz}
\usepackage{hyperref}
\usepackage{verbatim}
\usepackage{cleveref}
\usepackage{caption}
\usepackage{subcaption}
\usepackage[normalem]{ulem}
\usepackage{float}

\newtheorem{theorem}{Theorem}[section]
\newtheorem{proposition}[theorem]{Proposition}
\newtheorem{lemma}[theorem]{Lemma}

\newtheorem{corollary}[theorem]{Corollary}

\newtheorem{remark}[theorem]{Remark}

\newcommand{\R}{\mathbb{R}}
\newcommand{\Z}{\mathbb{Z}}
\newcommand{\N}{\mathbb{N}}
\newcommand{\E}{\mathbb{E}}
\newcommand{\Set}{\text{Set}}

\definecolor{lumicolor}{HTML}{FF6EFF}
\definecolor{cmagenta}{HTML}{CC79A7} 
\definecolor{corange}{HTML}{E69F00} 
\definecolor{cblue}{HTML}{0072B2} 

\newcommand{\OFG}{{\rm{OFG}}}

\title{The Origami flip graph of the $2\times n$ Miura-ori}
\author[Christensen, Hull, O'Neil, Pappano, Ter-Saakov, Yang]{Lumi Christensen}
\address{University of Pennsylvania}
\email{lumic@sas.upenn.edu}

\author[]{Thomas C. Hull}
\address{Franklin and Marshall College}
\email{thull1@fandm.edu}

\author[]{Emma O'Neil}
\address{Portland State University}
\email{emmao@pdx.edu}

\author[]{Valentina Pappano} 
\address{Rutgers University}
\email{valentina.pappano@gmail.com}

\author[]{Natalya Ter-Saakov}
\address{Rutgers University}
\email{nt399@rutgers.edu}

\author[]{Kacey Yang}
\address{University of California, Los Angeles}
\email{kaceyyang@g.ucla.edu}

\date{}

\begin{document}

\begin{abstract}
Given an \textit{origami crease pattern} $C=(V,E)$, a straight-line planar graph embedded in a region of $\R^2$, we assign each crease to be either a mountain crease (which bends convexly) or a valley crease (which bends concavely), creating a \textit{mountain-valley (MV) assignment} $\mu:E\to\{-1,1\}$. An MV assignment $\mu$ is \textit{locally valid} if the faces around each vertex in $C$ can be folded flat under $\mu$. In this paper, we investigate locally valid MV assignments of the \textit{Miura-ori}, $M_{m,n}$, an $m\times n$  parallelogram tessellation used in numerous engineering applications. The \textit{origami flip graph} $\OFG(C)$ of $C$ is a graph whose vertices are locally valid MV assignments of $C$, and two vertices are adjacent if they differ by a \textit{face flip}, an operation that swaps the MV-parity of every crease bordering a given face of $C$.  We enumerate the number of vertices and edges in $\OFG(M_{2,n})$ and prove several facts about the degrees of vertices in $\OFG(M_{2,n})$. By finding recurrence relations, we show that the number of vertices of degree $d$ and $2n-a$ (for $0\leq a$) are both described by polynomials of particular degrees. We then prove that the diameter of $\OFG(M_{2,n})$ is $\lceil \frac{n^2}{2}\rceil$ and find lower bounds on the diameter of $\OFG(M_{m,n})$ using techniques from $3$-coloring reconfiguration graphs.
\end{abstract}

\maketitle



\section{Introduction}\label{sec:intro}

The \textit{Miura-ori} or \textit{Miura fold} \cite{originalMiura} is an origami crease pattern with numerous engineering and architecture applications due in part to its ability to fold into a flat, compact shape in one continuous motion.  
For example, NASA has explored using origami crease patterns inspired by the Miura-ori to compactly fold solar arrays and then unfold them after being deployed in space \cite{Starshade}. The primary objective of this paper is to investigate the properties and structure of the Miura-ori \textit{origami flip graph}, a graph that encodes the relationship between valid folded configurations of the Miura-ori crease pattern (defined in Section~\ref{sec:bg}). Origami flip graphs were first introduced in \cite{tomTessellations}, and expanded in \cite{natashaOFG} with respect to single-vertex crease patterns. 
As we will see, the origami flip graph describes how folded configurations are connected via \textit{face flips}, an operation that has been used to study applications of the Miura-ori in metamaterials \cite{Silverberg1} and statistical mechanics \cite{Assis}.  Face flips are useful because they describe a local alteration of the crease pattern's configuration that can affect physical properties of the origami. For example, in \cite{Silverberg1} face flips are shown to be a simple way to mechanically tune the Miura-ori's resistance force under compression. Knowing more about the origami flip graph's structure, such as its connectivity and diameter, would give engineers a bound on the number of face flips needed to reprogram a Miura-ori's compression stiffness from one setting to any other. 
Although there has been significant study of the Miura-ori, there is very little literature on its origami flip graph. Currently, the only published result is that the origami flip graph of the Miura-ori is connected \cite{tomTessellations}. Thus, this paper represents the first deeper investigation into this area. 

While some of our results hold for the general $m\times n$ Miura-ori, which we denote $M_{m,n}$,
we concentrate on the origami flip graph of the Miura-ori with two rows and $n$ columns of parallelograms, $M_{2,n}$. Oftentimes, adding a new row to the Miura-ori would interfere with the results derived from the $2 \times n$ case, as the new row of vertices introduced would share diagonal creases with the previous row, thus removing a lot of the ``independence" of each vertex and making cases a lot more complex to handle. Hence it is not clear whether or not results from the $2 \times n$ case will generalize; we will try to specify our thoughts on this throughout. 
After establishing background information and formal definitions in Section \ref{sec:bg}, we derive a closed-form expression for the number of vertices and edges of the origami flip graph of $M_{2,n}$, denoted $\OFG(M_{2,n})$, in Section \ref{sec:verticesedges}. We then consider extending the MV assignments of the $2 \times n$ Miura-ori to a $2 \times (n+1)$ Miura-ori and use these extensions to prove several properties about the degree sequences of $\OFG(M_{2,n})$ in Section \ref{subsec:extensions}. Finally in Section \ref{subsec:diameter}, we calculate the diameter of $\OFG(M_{2,n})$ using a previously-known bijection between the Miura-ori and 3-colored grid graphs from \cite{Ginepro}.


\section{Background}\label{sec:bg}
\begin{figure}[!ht]
\centering
\begin{subfigure}[b]{0.4\textwidth}
    \centering
    \resizebox{\textwidth}{!}{%
\begin{circuitikz}
\tikzstyle{every node}=[font=\LARGE]
\draw  (12.5,9.75) -- (16.25,9.75) -- (15,7.25) -- (11.25,7.25) -- cycle;
\draw  (16.25,9.75) -- (20,9.75) -- (18.75,7.25) -- (15,7.25) -- cycle;
\draw  (20,9.75) -- (23.75,9.75) -- (22.5,7.25) -- (18.75,7.25) -- cycle;
\draw  (23.75,9.75) -- (27.5,9.75) -- (26.25,7.25) -- (22.5,7.25) -- cycle;
\draw  (12.5,4.75) -- (16.25,4.75) -- (15,2.25) -- (11.25,2.25) -- cycle;
\draw  (16.25,4.75) -- (20,4.75) -- (18.75,2.25) -- (15,2.25) -- cycle;
\draw  (20,4.75) -- (23.75,4.75) -- (22.5,2.25) -- (18.75,2.25) -- cycle;
\draw  (23.75,4.75) -- (27.5,4.75) -- (26.25,2.25) -- (22.5,2.25) -- cycle;
\draw  (11.25,7.25) -- (15,7.25) -- (16.25,4.75) -- (12.5,4.75) -- cycle;
\draw  (15,7.25) -- (18.75,7.25) -- (20,4.75) -- (16.25,4.75) -- cycle;
\draw  (18.75,7.25) -- (22.5,7.25) -- (23.75,4.75) -- (20,4.75) -- cycle;
\draw  (22.5,7.25) -- (26.25,7.25) -- (27.5,4.75) -- (23.75,4.75) -- cycle;
\draw  (11.25,2.25) -- (15,2.25) -- (16.25,-0.25) -- (12.5,-0.25) -- cycle;
\draw  (15,2.25) -- (18.75,2.25) -- (20,-0.25) -- (16.25,-0.25) -- cycle;
\draw  (18.75,2.25) -- (22.5,2.25) -- (23.75,-0.25) -- (20,-0.25) -- cycle;
\draw  (22.5,2.25) -- (26.25,2.25) -- (27.5,-0.25) -- (23.75,-0.25) -- cycle;
\end{circuitikz}
}%
\caption{A $4 \times 4$ Miura-ori crease pattern.}
\label{fig:miura}
\end{subfigure}
\hfill
\begin{subfigure}[b]{0.57\textwidth}
    \centering
\resizebox{\textwidth}{!}{
\begin{circuitikz}
\tikzstyle{every node}=[font=\LARGE]

\draw [line width=0.7pt, short] (7.5,6) -- (11.25,6);
\draw [line width=0.7pt, short] (7.5,6) -- (6.25,3.5);
\draw [line width=0.7pt, short] (11.25,6) -- (10,3.5)node[pos=0.5, fill=white]{$e_2$};
\draw [line width=0.7pt, short] (6.25,3.5) -- (10,3.5)node[pos=0.5, fill=white]{$e_0=e_1$};
\draw [line width=0.7pt, short] (6.25,3.5) -- (7.5,1);
\draw [line width=0.7pt, short] (7.5,1) -- (11.25,1);
\draw [line width=0.7pt, short] (10,3.5) -- (11.25,1)node[pos=0.5,fill=white]{$e_4$};
\draw [line width=0.7pt, short] (10,3.5) -- (13.75,3.5)node[pos=0.5, fill=white]{$e_3$};
\draw [line width=0.7pt, short] (11.25,6) -- (15,6);
\draw [line width=0.7pt, short] (11.25,1) -- (15,1);
\draw [line width=0.7pt, short] (15,6) -- (13.75,3.5)node[pos=0.5, fill=white]{$e_5$};
\draw [line width=0.7pt, short] (13.75,3.5) -- (17.5,3.5)node[pos=0.5, fill=white]{$e_6$};
\draw [line width=0.7pt, short] (13.75,3.5) -- (15,1)node[pos=0.5, fill=white]{$e_7$};
\draw [line width=0.7pt, short] (15,6) -- (18.75,6);
\draw [line width=0.7pt, short] (18.75,6) -- (17.5,3.5);
\draw [line width=0.7pt, short] (17.5,3.5) -- (18.75,1);
\draw [line width=0.7pt, short] (15,1) -- (18.75,1);
\draw  (19.75,3.5) circle (0cm);
\draw [dashed] (18.75,6) -- (23.75,6);
\draw [dashed] (17.5,3.5) -- (22.5,3.5);
\draw [dashed] (18.75,1) -- (23.75,1);
\draw [dashed] (22.5,6) -- (21.25,3.5);
\node [font=\LARGE, color={rgb,255:red,250; green,0; blue,0}] at (8.5,4.75) {$\alpha_{1,1}$};
\node [font=\LARGE, color={rgb,255:red,250; green,0; blue,0}] at (12.5,4.75) {$\alpha_{1,2}$};
\node [font=\LARGE, color={rgb,255:red,250; green,0; blue,0}] at (16.25,4.75) {$\alpha_{1,3}$};
\node [font=\LARGE, color={rgb,255:red,250; green,0; blue,0}] at (8.5,2.25) {$\alpha_{2,1}$};
\node [font=\LARGE, color={rgb,255:red,250; green,0; blue,0}] at (12.3,2.25) {$\alpha_{2,2}$};
\node [font=\LARGE, color={rgb,255:red,250; green,0; blue,0}] at (16.2,2.25) {$\alpha_{2,3}$};
\draw [dashed] (21.25,3.5) -- (22.5,1);
\node [font=\LARGE, color={rgb,255:red,0; green,0; blue,250}] at (9.6,4) {$x_1$};
\node [font=\LARGE, color={rgb,255:red,0; green,0; blue,250}] at (13.35,4) {$x_2$};
\end{circuitikz}}
\caption{Labelings of the $2\times n$ Miura-ori.}
\label{fig:edgelabel}
\end{subfigure}
    \caption{The Miura-ori crease pattern.}
\end{figure}
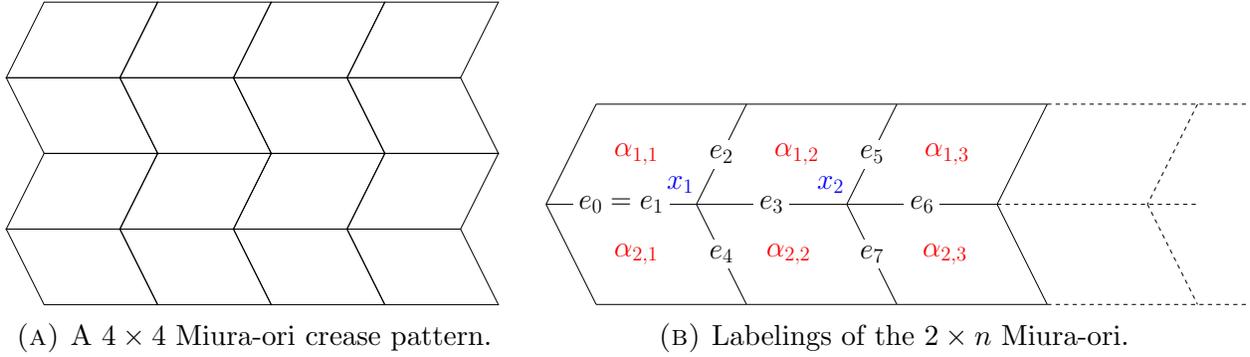

An \textit{origami crease pattern} is a straight line drawing of a planar graph $C=(V,E)$ on a region $A \subseteq\R^2$. The \textit{Miura-ori} \cite{originalMiura} is a crease pattern formed by crossing parallel lines with zig-zags, creating faces that are all congruent parallelograms (see Figure \ref{fig:miura}). We denote an $m$ by $n$ Miura-ori by $M_{m,n}$ and denote the face in the $i$th row, $j$th column as $\alpha_{i,j}$. The case $M_{2,n}$ will have a single row of vertices, which we call $x_1,...,x_{n-1}$ from left to right. For $k=1,\ldots, n-1$, we label the crease left of $x_k$ as $e_{3k-3}$, the crease above as $e_{3k-1}$, the crease to the right as $e_{3k}$, and the crease below as $e_{3k+1}$, as shown in Figure \ref{fig:edgelabel}. We additionally set $e_1=e_0$ to describe the leftmost crease. 

In this paper we will sometimes discuss the full $M_{m,n}$ Miura-ori, and in that case we denote $x_{i,j}$ to be the $i$th row, $j$th column vertex in the crease pattern and let $\{\alpha_{i,j},\alpha_{k,l}\}$ denote the crease line between the faces $\alpha_{i,j}$ and $\alpha_{k,l}$, if it exists.

A \textit{mountain-valley (MV) assignment} on a crease pattern $C=(V,E)$ is a function $\mu: E \to \{-1,1\}$, where $\mu(e)=-1$ means the crease $e$ is a \textit{valley} (concave) crease and $\mu(e)=1$ means  a \textit{mountain} (convex) crease. In this paper, we use solid lines to indicate mountain creases and dashed lines to indicate valleys. We say a crease pattern is \textit{locally flat-foldable} if there exists an MV assignment such that 
for each interior vertex $x$ in $C$, the faces of $C$ that are adjacent to $x$ can be folded flat without self-intersections (i.e., viewed as a crease pattern by itself, the vertex $x$ folds flat with its Ms and Vs). An MV assignment that achieves this is called \textit{locally valid}. 

In general, describing what MV assignments are flat-foldable and which are not is difficult \cite{boxpleat}. Luckily the single-vertex case is completely known:

\begin{theorem}[Kawasaki's Theorem \cite{tomMaekawa}]
    Let $C$ be a crease pattern with one interior vertex $x$ and sector angles between consecutive creases $\rho_1,\ldots,\rho_{2k}$. Then, $C$ is locally flat-foldable if and only if $$\rho_1-\rho_2+\rho_3-...-\rho_{2k}=0.$$
\end{theorem}

The mountains and valleys at a flat-foldable vertex must satisfy an invariant known as \textit{Maekawa's Theorem}:

\begin{theorem}[Maekawa's Theorem \cite{tomMaekawa,originalMaekawa}]
    Let $x$ be a vertex in a crease pattern with locally valid MV assignment $\mu$ and let $E$ be the set of creases incident to $x$. Then, 
    $$\sum_{e\in E} \mu(e)=\pm 2$$
\end{theorem}
This theorem states that around any vertex, the number of mountain creases and valley creases have to differ by 2. In the case where all sector angles between creases are equal, Maekawa's Theorem is necessary and sufficient to determine whether an MV assignment is flat-foldable at the vertex \cite{tomTessellations}. When the sector angles are not all equal, the smallest ones constrain the mountains and valleys as follows: 
\begin{theorem}[Big-Little-Big Lemma \cite{tomMaekawa}] \label{biglittlebig}
    Let $x$ be a vertex in a flat-foldable crease pattern with locally valid MV assignment $\mu$, and suppose that between the crease edges $e_i,...,e_{i+k+1}$ at $x$, the sector angles $\rho_i$ (between $e_i$ and $e_{i+1}$) satisfy the condition $\rho_{i+1}=...=\rho_{i+k}$, with $\rho_{i}>\rho_{i+1}$ and $\rho_{i+k+1}>\rho_i$. Then,
    $$\sum_{j=i+1}^{i+k+1}\mu(e_j)=\begin{cases}
        0, & \text{if k is odd}\\
        \pm 1, & \text{if k is even.}
    \end{cases}$$
\end{theorem}
That is, if we have an even number of small equal sector angles around a vertex that are bordered by larger angles, then among the edges that border these little sector angles, there needs to either be one more mountain crease than valley crease or vice versa. If we have an odd number of small sector angles, then there need to be equal numbers of mountain and valley creases. 

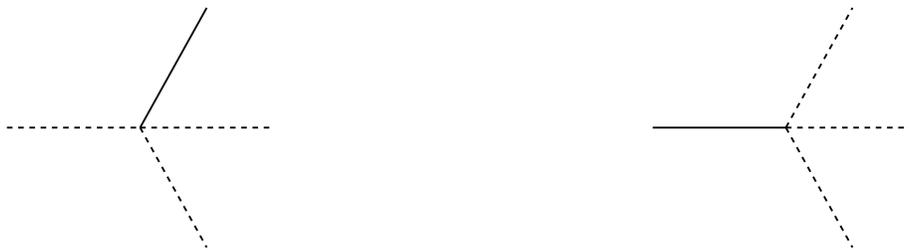
\begin{figure}
    \centering
    \begin{subfigure}[b]{0.48\textwidth}
        \centering
        \resizebox{0.45\textwidth}{!}{%
\begin{circuitikz}
\tikzstyle{every node}=[font=\LARGE]
\node  at (7.5,7) {(a)};
\draw [line width=1pt, short] (10,4.75) -- (11.25,7);
\draw [line width=1pt, dashed] (7.5,4.75) -- (10,4.75);
\draw [line width=1pt, dashed] (10,4.75) -- (12.5,4.75);
\draw [line width=1pt, dashed] (10,4.75) -- (11.25,2.5);
\end{circuitikz}}%
    \end{subfigure}
    \hfill
    \begin{subfigure}[b]{0.48\textwidth}
        \centering
        \resizebox{0.45\textwidth}{!}{%
        \begin{circuitikz}
            \tikzstyle{every node}=[font=\LARGE]
\node  at (13.5,7) {(b)};
\draw [line width=1pt, dashed] (16,4.75) -- (17.25,7);
\draw [line width=1pt, dashed] (16,4.75) -- (18.5,4.75);
\draw [line width=1pt, dashed] (16,4.75) -- (17.25,2.5);
\draw [line width=1pt, short] (13.5,4.75) -- (16,4.75);
        \end{circuitikz}}%
    \end{subfigure}
    \caption{Two MV assignments of a Miura-ori vertex: (a) a valid and (b) an invalid MV assignment. The three ``toes" of a vertex cannot all be the same parity. }
    \label{valid MV}
    \end{figure}

We use these theorems to determine the conditions a locally valid MV assignment on the Miura-ori must satisfy. Each vertex of $M_{m,n}$ has 4 sector angles (2 obtuse, 2 acute) bordered by 4 creases. By Kawasaki's Theorem, we know that the placement of these two obtuse and two acute angles ensure that the Miura-ori is locally flat-foldable. By Maekawa's Theorem, we know that the four creases around each vertex either have to be 3 valley creases and 1 mountain crease, or 3 mountain creases and 1 valley crease. By the Big-Little-Big Lemma, we know that the minority crease (the crease with the different MV assignment) cannot be the one bordered by two obtuse angles. This is because it would cause the three creases (the three ``toes" of the vertex) that border the two acute angles to all be the same parity, which violates the Big-Little-Big Lemma (see Figure \ref{valid MV}). 

Given a face $\alpha$ of a flat-foldable crease pattern and a MV assignment $\mu$, a  \textit{face flip of $\alpha$ under $\mu$} is a new MV assignment $\mu_\alpha$ where we ``flip" the edges bordering  $\alpha$ to obtain $\mu_\alpha$.  That is,  $\mu_\alpha(e)=-\mu(e)$, if $e$ borders $\alpha$, and $\mu_\alpha(e)=\mu(e)$, otherwise.
If $\mu$ and $\mu_\alpha$ are both locally valid MV assignments, then we consider $\alpha$ to be a \textit{flippable face under $\mu$}. 

\begin{figure}
\centering
\resizebox{0.6\textwidth}{!}{%
\begin{circuitikz}
\tikzstyle{every node}=[font=\Huge]
\draw [line width=1.5pt, short] (6.25,20) -- (5,17.5);
\draw [line width=1.5pt, short] (5,17.5) -- (6.25,15);
\draw [line width=1.5pt, short] (6.25,20) -- (11.25,20);
\draw [line width=1.5pt, short] (11.25,20) -- (10,17.5);
\draw [line width=1.5pt, short] (10,17.5) -- (11.25,15);
\draw [line width=1.5pt, short] (11.25,15) -- (6.25,15);
\draw [line width=1.5pt, short] (15,20) -- (13.75,17.5);
\draw [line width=1.5pt, short] (13.75,17.5) -- (15,15);
\draw [line width=1.5pt, short] (15,15) -- (20,15);
\draw [line width=1.5pt, short] (20,15) -- (18.75,17.5);
\draw [line width=1.5pt, short] (18.75,17.5) -- (20,20);
\draw [line width=1.5pt, short] (20,20) -- (15,20);
\draw [line width=1.5pt, short] (8.75,20) -- (7.5,17.5);
\draw [line width=1.5pt, short] (7.5,17.5) -- (10,17.5);
\draw [line width=1.5pt, short] (7.5,17.5) -- (5,17.5);
\draw [line width=1.5pt, dashed] (7.5,17.5) -- (8.75,15);
\draw [line width=1.5pt, dashed] (18.75,17.5) -- (16.25,17.5);
\draw [line width=1.5pt, dashed] (16.25,17.5) -- (17.5,15);
\draw [line width=1.5pt, dashed] (16.25,17.5) -- (13.75,17.5);
\draw [line width=1.5pt, short] (17.5,20) -- (16.25,17.5);
\node [font=\Huge, color={rgb,255:red,0; green,200; blue,50}] at (9.5,16.25) {$\checkmark$};
\node [font=\Huge, color={rgb,255:red,0; green,200; blue,50}] at (6.75,18.75) {$\checkmark$};
\node [font=\Huge, color={rgb,255:red,0; green,200; blue,50}] at (18.25,18.75) {$\checkmark$};
\node [font=\Huge, color={rgb,255:red,0; green,200; blue,50}] at (15.5,16.25) {$\checkmark$};
\draw [line width=1.5pt, short] (6.25,13.75) -- (5,11.25);
\draw [line width=1.5pt, short] (5,11.25) -- (6.25,8.75);
\draw [line width=1.5pt, short] (6.25,13.75) -- (11.25,13.75);
\draw [line width=1.5pt, short] (11.25,13.75) -- (10,11.25);
\draw [line width=1.5pt, short] (10,11.25) -- (11.25,8.75);
\draw [line width=1.5pt, short] (11.25,8.75) -- (6.25,8.75);
\draw [line width=1.5pt, short] (15,13.75) -- (13.75,11.25);
\draw [line width=1.5pt, short] (13.75,11.25) -- (15,8.75);
\draw [line width=1.5pt, short] (15,8.75) -- (20,8.75);
\draw [line width=1.5pt, short] (20,8.75) -- (18.75,11.25);
\draw [line width=1.5pt, short] (18.75,11.25) -- (20,13.75);
\draw [line width=1.5pt, short] (20,13.75) -- (15,13.75);
\node [font=\Huge, color={rgb,255:red,0; green,200; blue,50}] at (9.5,10) {$\checkmark$};
\node [font=\Huge, color={rgb,255:red,0; green,200; blue,50}] at (6.75,12.5) {$\checkmark$};
\node [font=\Huge, color={rgb,255:red,0; green,200; blue,50}] at (18.25,12.5) {$\checkmark$};
\node [font=\Huge, color={rgb,255:red,0; green,200; blue,50}] at (15.5,10) {$\checkmark$};
\draw [line width=1.5pt, dashed] (8.75,13.75) -- (7.5,11.25);
\draw [line width=1.5pt, dashed] (7.5,11.25) -- (5,11.25);
\draw [line width=1.5pt, dashed] (10,11.25) -- (7.5,11.25);
\draw [line width=1.5pt, dashed] (17.5,13.75) -- (16.25,11.25);
\draw [line width=1.5pt, short] (13.75,11.25) -- (18.75,11.25);
\draw [line width=1.5pt, short] (17.5,8.75) -- (16.25,11.25);
\draw [line width=1.5pt, short] (8.75,8.75) -- (7.5,11.25);
\draw [line width=1.5pt, short] (23.75,20) -- (22.5,17.5);
\draw [line width=1.5pt, short] (22.5,17.5) -- (23.75,15);
\draw [line width=1.5pt, short] (23.75,20) -- (28.75,20);
\draw [line width=1.5pt, short] (28.75,20) -- (27.5,17.5);
\draw [line width=1.5pt, short] (27.5,17.5) -- (28.75,15);
\draw [line width=1.5pt, short] (28.75,15) -- (23.75,15);
\draw [line width=1.5pt, short] (23.75,13.75) -- (22.5,11.25);
\draw [line width=1.5pt, short] (22.5,11.25) -- (23.75,8.75);
\draw [line width=1.5pt, short] (23.75,13.75) -- (28.75,13.75);
\draw [line width=1.5pt, short] (28.75,13.75) -- (27.5,11.25);
\draw [line width=1.5pt, short] (27.5,11.25) -- (28.75,8.75);
\draw [line width=1.5pt, short] (28.75,8.75) -- (23.75,8.75);
\draw [line width=1.5pt, dashed] (26.25,13.75) -- (25,11.25);
\draw [line width=1.5pt, dashed] (25,11.25) -- (22.5,11.25);
\draw [line width=1.5pt, short] (27.5,11.25) -- (25,11.25);
\draw [line width=1.5pt, dashed] (26.25,8.75) -- (25,11.25);
\draw [line width=1.5pt, short] (26.25,20) -- (25,17.5);
\draw [line width=1.5pt, short] (25,17.5) -- (22.5,17.5);
\draw [line width=1.5pt, dashed] (27.5,17.5) -- (25,17.5);
\draw [line width=1.5pt, short] (26.25,15) -- (25,17.5);
\node [font=\Huge, color={rgb,255:red,0; green,200; blue,50}] at (27,10) {$\checkmark$};
\node [font=\Huge, color={rgb,255:red,0; green,200; blue,50}] at (24.25,12.5) {$\checkmark$};
\node [font=\Huge, color={rgb,255:red,0; green,200; blue,50}] at (24.25,10) {$\checkmark$};
\node [font=\Huge, color={rgb,255:red,0; green,200; blue,50}] at (27,12.5) {$\checkmark$};
\node [font=\Huge, color={rgb,255:red,0; green,200; blue,50}] at (27,18.75) {$\checkmark$};
\node [font=\Huge, color={rgb,255:red,0; green,200; blue,50}] at (24.25,16.25) {$\checkmark$};
\node [font=\Huge, color={rgb,255:red,0; green,200; blue,50}] at (24.25,18.75) {$\checkmark$};
\node [font=\Huge, color={rgb,255:red,0; green,200; blue,50}] at (27,16.25) {$\checkmark$};
\end{circuitikz}
}%
\caption{MV assignments of a vertex with flippable faces indicated by check marks. Notice that around a single vertex, only diagonal faces can be made unflippable.}
\label{flippables}
\end{figure}
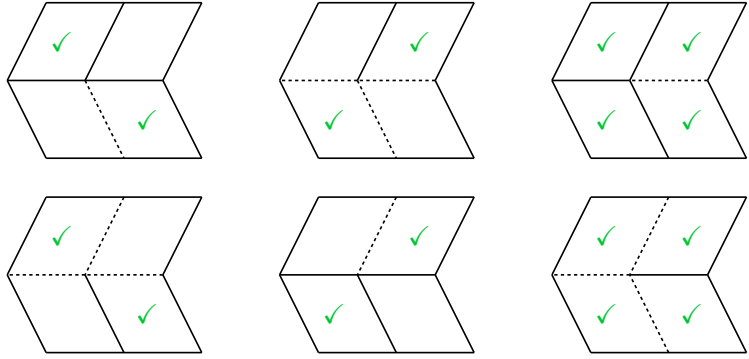

Figure~\ref{flippables} shows how the MV assignment $\mu$ at a Miura-ori vertex will affect a face's flippable nature. Let $e_j$ and $e_k$ be the two creases at this vertex that border a face $\alpha$. First note that if $\mu$ satisfies Maekawa's Theorem at the vertex, then $\mu_\alpha$ will also satisfy Maekawa at the vertex. Thus, to determine if $\alpha$ is flippable, we need only check that the Big-Little-Big Lemma is satisfied after flipping. If $e_j, e_k$ form an acute angle, then in order for $\alpha$ to be flippable we must have $\mu(e_j)\neq\mu(e_k)$, since if they were the same, flipping $\alpha$ would result in the three creases around the acute angles of the vertex being the same M/V parity.  Similarly, if $e_j,e_k$ meet at an obtuse angle, then we need $\mu(e_j)=\mu(e_k)$ to make $\alpha$ flippable. 

Let $C$ be a crease pattern. The \textit{origami flip graph of $C$}, denoted $\OFG(C)$, is a graph whose vertices are the locally valid MV assignments of $C$, and where two vertices are adjacent if and only if they differ by one face flip. Essentially, this graph provides us with an intuition on how ``close" certain MV assignments are to each other. An example of an origami flip graph is in Figure \ref{fig:2x2OFG}.

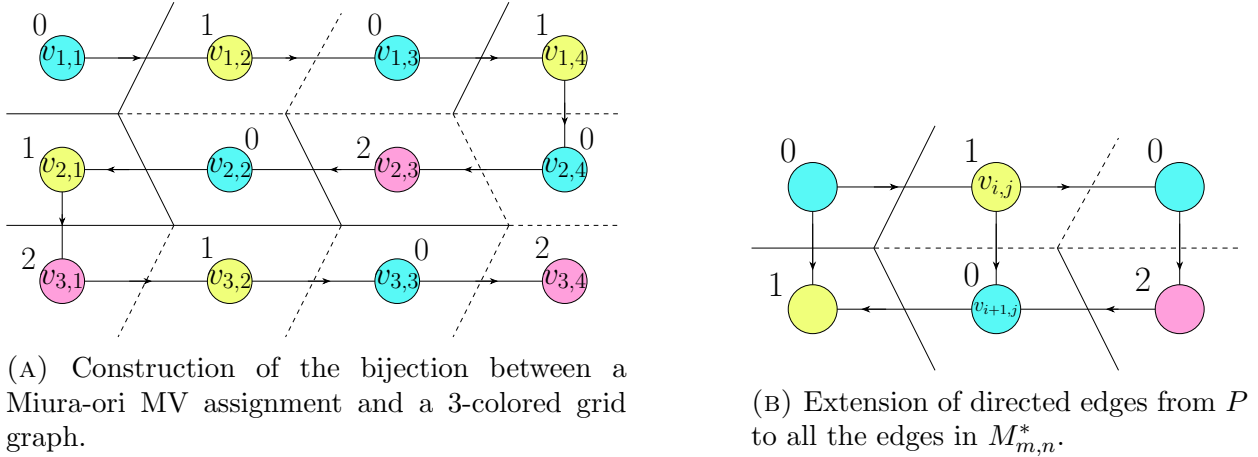
\begin{figure}
\centering
\begin{subfigure}[b]{0.5\textwidth}
    \centering
\resizebox{1\textwidth}{!}{%
\begin{circuitikz}
\tikzstyle{every node}=[font=\LARGE]
\draw [short] (11.25,9.75) -- (10,7.25);
\draw [short] (10,7.25) -- (11.25,4.75);
\draw [dashed] (10,7.25) -- (13.75,7.25);
\draw [dashed] (13.75,7.25) -- (15,9.5);
\draw [dashed] (13.75,7.25) -- (17.5,7.25);
\draw [short] (13.75,7.25) -- (15,4.75);
\draw [short] (18.75,9.75) -- (17.5,7.25);
\draw [dashed] (17.5,7.25) -- (21.25,7.25);
\draw [dashed] (17.5,7.25) -- (18.75,4.75);
\draw [dashed] (18.75,4.75) -- (21.25,4.75);
\draw [dashed] (18.75,4.75) -- (17.5,2.25);
\draw [short] (15,4.75) -- (18.75,4.75);
\draw [short] (11.25,4.75) -- (15,4.75);
\draw [dashed] (15,4.75) -- (13.75,2.25);
\draw [short] (7.5,7.25) -- (10,7.25);
\draw [ fill={rgb,255:red,88; green,249; blue,246} ] (8.75,8.5) circle (0.5cm) node {\LARGE $v_{1,1}$} ;
\draw [ fill={rgb,255:red,88; green,249; blue,246} ] (16.25,8.5) circle (0.5cm) node {\LARGE $v_{1,3}$} ;
\draw [ fill={rgb,255:red,88; green,249; blue,246} ] (20,6) circle (0.5cm) node {\LARGE $v_{2,4}$} ;
\draw [ fill={rgb,255:red,88; green,249; blue,246} ] (12.5,6) circle (0.5cm) node {\LARGE $v_{2,2}$} ;
\draw [ fill={rgb,255:red,255; green,159; blue,220} ] (8.75,3.5) circle (0.5cm) node {\LARGE $v_{3,1}$} ;
\draw [ fill={rgb,255:red,88; green,249; blue,246} ] (16.25,3.5) circle (0.5cm) node {\LARGE $v_{3,3}$} ;
\draw [ fill={rgb,255:red,255; green,159; blue,220} ] (16.25,6) circle (0.5cm) node {\LARGE $v_{2,3}$} ;
\draw [ fill={rgb,255:red,255; green,159; blue,220} ] (20,3.5) circle (0.5cm) node {\LARGE $v_{3,4}$} ;
\draw [ fill={rgb,255:red,240; green,255; blue,122} ] (12.5,8.5) circle (0.5cm) node {\LARGE $v_{1,2}$} ;
\draw [ fill={rgb,255:red,240; green,255; blue,122} ] (20,8.5) circle (0.5cm) node {\LARGE $v_{1,4}$} ;
\draw [ fill={rgb,255:red,240; green,255; blue,122} ] (8.75,6) circle (0.5cm) node {\LARGE $v_{2,1}$} ;
\draw [ fill={rgb,255:red,240; green,255; blue,122} ] (12.5,3.5) circle (0.5cm) node {\LARGE $v_{3,2}$} ;
\draw [short] (7.5,4.75) -- (11.25,4.75);
\draw [dashed] (11.25,4.75) -- (10,2.25);
\draw [short] (9.25,8.5) -- (12,8.5);
\draw [short] (13,8.5) -- (15.75,8.5);
\draw [short] (16.75,8.5) -- (19.5,8.5);
\draw [short] (20,8) -- (20,6.5);
\draw [short] (19.5,6) -- (16.75,6);
\draw [short] (15.75,6) -- (13,6);
\draw [short] (12,6) -- (9.25,6);
\draw [short] (8.75,5.5) -- (8.75,4);
\draw [short] (9.25,3.5) -- (12,3.5);
\draw [short] (13,3.5) -- (15.75,3.5);
\draw [short] (16.75,3.5) -- (19.5,3.5);
\draw [->, >=Stealth] (10,8.5) -- (10.5,8.5);
\draw [->, >=Stealth] (13.75,8.5) -- (14,8.5);
\draw [->, >=Stealth] (18,8.5) -- (18.5,8.5);
\draw [->, >=Stealth] (20,7.5) -- (20,7);
\draw [->, >=Stealth] (18.25,6) -- (17.75,6);
\draw [->, >=Stealth] (8.75,5) -- (8.75,4.75);
\draw [->, >=Stealth] (10.25,3.5) -- (10.75,3.5);
\draw [->, >=Stealth] (14.25,3.5) -- (14.75,3.5);
\draw [->, >=Stealth] (18,3.5) -- (18.5,3.5);
\draw [->, >=Stealth] (10.25,6) -- (9.75,6);
\draw [->, >=Stealth] (15,6) -- (14.75,6);
\node [font=\LARGE] at (8.25,9.25) {$0$};
\node [font=\LARGE] at (15.75,9.25) {$0$};
\node [font=\LARGE] at (20.5,6.7) {$0$};
\node [font=\LARGE] at (13,6.7) {$0$};
\node [font=\LARGE] at (16.8,4.2) {$0$};
\node [font=\LARGE] at (12,9.25) {$1$};
\node [font=\LARGE] at (19.5,9.25) {$1$};
\node [font=\LARGE] at (8,6.5) {$1$};
\node [font=\LARGE] at (12,4.25) {$1$};
\node [font=\LARGE] at (15.5,6.5) {$2$};
\node [font=\LARGE] at (19.5,4.25) {$2$};
\node [font=\LARGE] at (8,4) {$2$};
\end{circuitikz}
}%
\caption{Construction of the bijection between a Miura-ori MV assignment and a 3-colored grid graph.}
\label{fig:3colorbiject}
\end{subfigure}
\hfill
\begin{subfigure}[b]{0.40\textwidth}
\centering
\resizebox{1\textwidth}{!}{%
\begin{circuitikz}
\tikzstyle{every node}=[font=\LARGE]

\draw [short] (11.25,9.75) -- (10,7.25);
\draw [short] (10,7.25) -- (11.25,4.75);
\draw [dashed] (10,7.25) -- (13.75,7.25);
\draw [dashed] (13.75,7.25) -- (15,9.5);
\draw [dashed] (13.75,7.25) -- (17.5,7.25);
\draw [short] (13.75,7.25) -- (15,4.75);
\draw [short] (7.5,7.25) -- (10,7.25);
\draw [ fill={rgb,255:red,88; green,249; blue,246} ] (8.75,8.5) circle (0.5cm)  ;
\draw [ fill={rgb,255:red,88; green,249; blue,246} ] (16.25,8.5) circle (0.5cm) ;
\draw [ fill={rgb,255:red,88; green,249; blue,246} ] (12.5,6) circle (0.5cm) node {\normalsize $v_{i+1,j}$} ;
\draw [ fill={rgb,255:red,255; green,159; blue,220} ] (16.25,6) circle (0.5cm) ;
\draw [ fill={rgb,255:red,240; green,255; blue,122} ] (12.5,8.5) circle (0.5cm) node {\Large $v_{i,j}$} ;
\draw [ fill={rgb,255:red,240; green,255; blue,122} ] (8.75,6) circle (0.5cm) ;
\draw [short] (9.25,8.5) -- (12,8.5);
\draw [short] (13,8.5) -- (15.75,8.5);
\draw [short] (15.75,6) -- (13,6);
\draw [short] (12,6) -- (9.25,6);
\draw [->, >=Stealth] (10,8.5) -- (10.5,8.5);
\draw [->, >=Stealth] (13.75,8.5) -- (14,8.5);
\draw [->, >=Stealth] (10.25,6) -- (9.75,6);
\draw [->, >=Stealth] (15,6) -- (14.75,6);
\node [font=\LARGE] at (8.25,9.25) {$0$};
\node [font=\LARGE] at (15.75,9.25) {$0$};
\node [font=\LARGE] at (12,6.7) {$0$};
\node [font=\LARGE] at (12,9.25) {$1$};
\node [font=\LARGE] at (8,6.5) {$1$};
\node [font=\LARGE] at (15.5,6.6) {$2$};
\draw [short] (16.25,8) -- (16.25,6.5);
\draw [short] (12.5,8) -- (12.5,6.5);
\draw [short] (8.75,8) -- (8.75,6.5);
\draw [->, >=Stealth] (8.75,7.75) -- (8.75,6.8);
\draw [->, >=Stealth] (12.5,7.75) -- (12.5,6.8);
\draw [->, >=Stealth] (16.25,7.75) -- (16.25,6.8);
\end{circuitikz}
}%
\caption{Extension of directed edges from $P$ to all the edges in $M^*_{m,n}$.}
\label{fig:3colorbiject2}
\end{subfigure}
\caption{Illustrating the Miura-ori MV assignment bijection with 3-colorings.}
\end{figure}

We next describe a bijection, from \cite{Ginepro}, between the set of all locally-valid MV assignments of $M_{m,n}$ and the set of all proper 3-vertex colorings of the $m\times n$ ``dual" grid graph $M^*_{m,n}$ with one vertex precolored. 
Position $M_{m,n}$ with the top left corner being obtuse, and overlay it with $M^*_{m,n}$ so that each vertex corresponds to a face of $M_{m,n}$. Our colors for $M^*_{m,n}$ will be the elements of $\Z_3$. We follow a directed path $P$ in $M^*_{m,n}$ along a boustrophedon pattern: Start at the top left with $v_{1,1}$ (corresponding to $\alpha_{1,1}$ in $M_{m,n}$) then travel to the right through all the vertices in the top row of $M^*_{m,n}$, then move down to the second row from $v_{1,n}$ to $v_{2,n}$ and then travel left through the vertices in the second row, then move down to the third row, travel right through all the vertices in the third row, and so on until all vertices are included in the path. We color the vertices as follows: Color the top left vertex $v_{1,1}$  $0\mod 3$. Then we color the vertices in $P$ consecutively by adding $\mu(e_k)$ to the previous vertex's color, where $e_k$ is the crease of $M_{m,n}$ crossed by the edge of $P$ (see Figure \ref{fig:3colorbiject}). More formally, our coloring $\gamma: V(M^*_{m,n})\to \Z_3$ is as follows: 
$$\begin{cases}
    \gamma(v_{1,1})=0\\
    \gamma(v_{i,j+1})=\gamma(v_{i,j})+\mu(e_k)\mod 3, & i \text{ odd}, 1\leq j\leq n-1\\
    \gamma(v_{i,j-1})=\gamma(v_{i,j})+\mu(e_k)\mod 3, & i \text{ even}, 2\leq j\leq n\\
\gamma(v_{i+1,n})=\gamma(v_{i,n})+\mu(e_k)\mod 3, & i\text{ odd}\\
\gamma(v_{i+1,1})=\gamma(v_{i,1})+\mu(e_k)\mod 3, & i\text{ even}
\end{cases}$$ 
It can be shown that this process constructs a proper $3$-colored grid graph for every locally valid MV assignment, and that every 3-coloring has a corresponding locally valid MV assignment \cite{Ginepro, tomOrigametry}. 

We extend our directed edges from $P$  
to direct all the edges in $M^*_{m,n}$ by orienting all the edges not in $P$ downward.  These new directed edges retain the property of crossing over creases leading to increases or decreases in color, as proved in \cite{Ginepro}. 
That is, if the directed edge from $v_{i,j}$ to $v_{i+1,j}$  
crosses over crease $e_\ell$, the following is still true (See Figure \ref{fig:3colorbiject2}):
$$\gamma(v_{i+1,j})=\gamma(v_{i,j})+\mu(e_\ell)\mod 3.$$

We fix the first vertex color arbitrarily. Thus, instead of thinking of the bijection as between MV assignments and three colorings of grid graphs with the first vertex color fixed, we may consider it instead as a bijection between MV assignments and 3-coloring equivalence classes, where any two elements in each equivalence class only differ by a color permutation. I.e., if $\gamma_1$ and $\gamma_2$ are two different colorings with $\gamma_1(v_{1,1})=0$, $\gamma_2(v_{1,1})\neq 0$, and 
$$\begin{cases}
    \gamma_2(v_{i+1,j})-\gamma_2(v_{i,j})=\gamma_1(v_{i+1,j})-\gamma_1(v_{i,j}),& 1\leq i\leq m-1, 1\leq j\leq n\\
    \gamma_2(v_{i,j+1})-\gamma_2(v_{i,j})=\gamma_1(v_{i,j+1})-\gamma_1(v_{i,j}),& 1\leq i\leq m, 1\leq j\leq n-1, 
\end{cases}$$ 
then they would be in the same equivalence class; $\mu$ would correspond to both 3-colorings. 

By this bijection, flipping a face $\alpha_{i,j}$ of a Miura-ori is equivalent to changing the color of the corresponding vertex $v_{i,j}$ in the $3$-colored grid graph, and a face of $M_{m,n}$ is flippable if and only if the color of the corresponding vertex of $M^*_{m,n}$ can be changed while remaining a proper $3$-coloring (see \cite[Lemma 4.1]{tomTessellations}).


\section{Number of vertices and edges in \texorpdfstring{$\OFG(M_{2,n})$}{OFG(M2,n)}}\label{sec:verticesedges}

\begin{figure}
\centering
\resizebox{0.8\textwidth}{!}{
\begin{circuitikz}
\tikzstyle{every node}=[font=\large]
\draw [line width=2pt, short] (17.5,20) -- (16.25,17.5);
\draw [line width=2pt, short] (16.25,17.5) -- (17.5,15);
\draw [line width=2pt, short] (17.5,15) -- (20,15);
\draw [line width=2pt, short] (20,15) -- (22.5,15);
\draw [line width=2pt, short] (22.5,15) -- (21.25,17.5);
\draw [line width=2pt, short] (21.25,17.5) -- (22.5,20);
\draw [line width=2pt, short] (22.5,20) -- (17.5,20);
\draw [line width=2pt, short] (20,20) -- (18.75,17.5);
\draw [line width=2pt, short] (18.75,17.5) -- (16.25,17.5);
\draw [short] (18.75,17.5) -- (20,15);
\draw [line width=2pt, dashed] (18.75,17.5) -- (21.25,17.5);
\node [font=\Huge, scale=1.5, color={rgb,255:red,64; green,200; blue,64}] at (18.1,18.75) {$\alpha_{1,1}$};
\node [font=\Huge, scale=1.5, color={rgb,255:red,180; green,180; blue,64}] at (20.6,18.75) {$\alpha_{1,2}$};
\node [font=\Huge, scale=1.5, color={rgb,255:red,200; green,64; blue,64}] at (20.5,16.25) {$\alpha_{2,2}$};
\node [font=\Huge, scale=1.5, color={rgb,255:red,64; green,64; blue,200}] at (18,16.25) {$\alpha_{2,1}$};
\draw [line width=2pt, short] (22.5,12.5) -- (21.25,10);
\draw [line width=2pt, short] (21.25,10) -- (22.5,7.5);
\draw [line width=2pt, short] (22.5,7.5) -- (27.5,7.5);
\draw [line width=2pt, short] (27.5,7.5) -- (26.25,10);
\draw [line width=2pt, short] (26.25,10) -- (27.5,12.5);
\draw [line width=2pt, short] (27.5,12.5) -- (22.5,12.5);
\draw [line width=2pt, short] (17.5,12.5) -- (12.5,12.5);
\draw [line width=2pt, short] (17.5,12.5) -- (16.25,10);
\draw [line width=2pt, short] (16.25,10) -- (17.5,7.5);
\draw [line width=2pt, short] (12.5,12.5) -- (11.25,10);
\draw [line width=2pt, short] (11.25,10) -- (12.5,7.5);
\draw [line width=2pt, short] (12.5,7.5) -- (17.5,7.5);
\draw [line width=2pt, short] (7.5,12.5) -- (2.5,12.5);
\draw [line width=2pt, short] (7.5,12.5) -- (6.25,10);
\draw [line width=2pt, short] (6.25,10) -- (7.5,7.5);
\draw [line width=2pt, short] (2.5,12.5) -- (1.25,10);
\draw [line width=2pt, short] (1.25,10) -- (2.5,7.5);
\draw [line width=2pt, short] (2.5,7.5) -- (7.5,7.5);
\draw [line width=2pt, short] (37.5,12.5) -- (32.5,12.5);
\draw [line width=2pt, short] (37.5,12.5) -- (36.25,10);
\draw [line width=2pt, short] (36.25,10) -- (37.5,7.5);
\draw [line width=2pt, short] (32.5,12.5) -- (31.25,10);
\draw [line width=2pt, short] (31.25,10) -- (32.5,7.5);
\draw [line width=2pt, short] (32.5,7.5) -- (37.5,7.5);
\draw [ color={rgb,255:red,200; green,64; blue,64}, line width=4pt, short] (3.75,10) -- (5,7.5);
\draw [ color={rgb,255:red,200; green,64; blue,64}, line width=6pt, dashed] (3.75,10) -- (6.25,10);
\draw [ color={rgb,255:red,64; green,200; blue,64}, line width=6pt, dashed] (3.75,10) -- (5,12.5);
\draw [ color={rgb,255:red,64; green,200; blue,64}, line width=6pt, dashed] (3.75,10) -- (1.25,10);
\draw [ color={rgb,255:red,64; green,200; blue,64}, line width=4pt, short] (15,12.5) -- (13.75,10);
\draw [ color={rgb,255:red,64; green,200; blue,64}, line width=4pt, short] (13.75,10) -- (11.25,10);
\draw [ color={rgb,255:red,200; green,64; blue,64}, line width=4pt, short] (13.75,10) -- (16.25,10);
\draw [ color={rgb,255:red,200; green,64; blue,64}, line width=6pt, dashed] (13.75,10) -- (15,7.5);
\draw [line width=2pt, short] (15.5,17.5) .. controls (9.5,18) and (5.5,15.75) .. (5,13)node[pos=0.5, fill=white, font=\Huge]{$\alpha_{1,1}$};
\draw [line width=2pt, short] (16.25,16.25) -- (15,13)node[pos=0.5, fill=white, font=\Huge]{$\alpha_{2,2}$};
\draw [line width=2pt, short] (22.5,16.25) -- (25,13)node[pos=0.5, fill=white, font=\Huge]{$\alpha_{1,2}$};
\draw [line width=2pt, short] (23.25,17.5) .. controls (29.5,18) and (34.25,17) .. (35,13)node[pos=0.5, fill=white, font=\Huge]{$\alpha_{2,1}$};
\draw [ color={rgb,255:red,64; green,64; blue,200}, line width=4pt, short] (21.25,10) -- (23.75,10);
\draw [ color={rgb,255:red,64; green,64; blue,200}, line width=4pt, short] (23.75,10) -- (25,7.5);
\draw [ color={rgb,255:red,180; green,180; blue,64}, line width=4pt, short] (23.75,10) -- (26.25,10);
\draw [ color={rgb,255:red,180; green,180; blue,64}, line width=6pt, dashed] (23.75,10) -- (25,12.5);
\draw [ color={rgb,255:red,180; green,180; blue,64}, line width=6pt, dashed] (36.25,10) -- (33.75,10);
\draw [ color={rgb,255:red,180; green,180; blue,64}, line width=4pt, short] (33.75,10) -- (35,12.5);
\draw [ color={rgb,255:red,64; green,64; blue,200}, line width=6pt, dashed] (33.75,10) -- (31.25,10);
\draw [ color={rgb,255:red,64; green,64; blue,200}, line width=6pt, dashed] (33.75,10) -- (35,7.5);
\draw [line width=2pt, short] (18.75,17.5) -- (20,15);
\draw [line width=2pt, short] (17.5,5) -- (16.25,2.5);
\draw [line width=2pt, short] (16.25,2.5) -- (17.5,0);
\draw [line width=2pt, short] (17.5,0) -- (20,0);
\draw [line width=2pt, short] (20,0) -- (22.5,0);
\draw [line width=2pt, short] (22.5,0) -- (21.25,2.5);
\draw [line width=2pt, short] (21.25,2.5) -- (22.5,5);
\draw [line width=2pt, short] (22.5,5) -- (17.5,5);
\node [font=\huge, scale=1.5, color={rgb,255:red,64; green,200; blue,64}] at (18.1,3.75) {$\alpha_{1,1}$};
\node [font=\huge, scale=1.5, color={rgb,255:red,180; green,180; blue,64}] at (20.6,3.75) {$\alpha_{1,2}$};
\node [font=\huge, scale=1.5, color={rgb,255:red,200; green,64; blue,64}] at (20.5,1.25) {$\alpha_{2,2}$};
\node [font=\huge, scale=1.5, color={rgb,255:red,64; green,64; blue,200}] at (18,1.25) {$\alpha_{2,1}$};
\draw [line width=2pt, dashed] (20,5) -- (18.75,2.5);
\draw [line width=2pt, dashed] (18.75,2.5) -- (16.25,2.5);
\draw [line width=2pt, dashed] (18.75,2.5) -- (20,0);
\draw [line width=2pt, short] (18.75,2.5) -- (21.25,2.5);
\draw [line width=2pt, short] (15,7) -- (16.25,3.75)node[pos=0.5, fill=white, font=\Huge]{$\alpha_{1,1}$};
\draw [line width=2pt, short] (5,7) .. controls (5.5,3) and (9.5,2) .. (15.75,2.5)node[pos=0.5, fill=white, font=\Huge]{$\alpha_{2,2}$};
\draw [line width=2pt, short] (25,6.75) -- (22.5,3.75)node[pos=0.5, fill=white, font=\Huge]{$\alpha_{2,1}$};
\draw [line width=2pt, short] (35,7) .. controls (34.25,3) and (30,2) .. (23.25,2.5)node[pos=0.5, fill=white, font=\Huge]{$\alpha_{1,2}$};
\end{circuitikz}
}
\caption{Origami flip graph of the $2\times 2$ Miura-ori, $\OFG(M_{2,2})$.}
\label{fig:2x2OFG}
\end{figure}
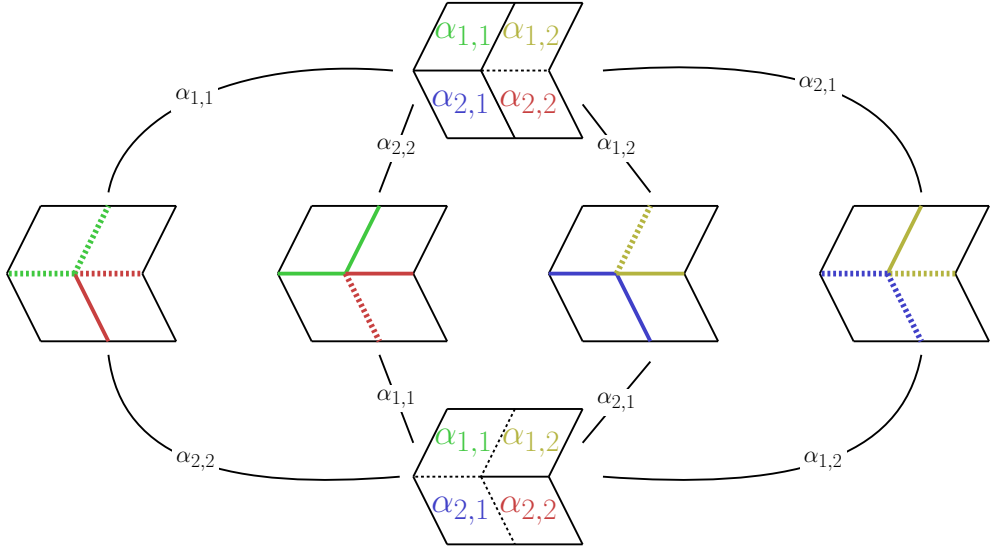

We now examine the structure of the origami flip graph of the $2\times n$ Miura-ori; Figure \ref{fig:2x2OFG} shows $\OFG(M_{2,2})$. We begin by enumerating the number of vertices and edges in $\OFG(M_{2,n})$. The number of vertices has already been shown in \cite{Ginepro}, but we repeat the result here for completeness. 

\begin{theorem}\label{thm:2xnvertices}
    The number of vertices of $\OFG(M_{2,n})$ is $2(3^{n-1})$. 
\end{theorem}
\begin{proof}
    We must count the number of valid MV assignments of $M_{2,n}$.
    Using the notation of Figure~\ref{fig:edgelabel}, crease $e_1$ has two choices for its M/V parity. 
    Then for the creases $e_2,e_3,e_4$, the Big-Little-Big and Maekawa's Theorems tell us that two of the creases have to match the parity of $e_1$ and the third has to be the opposite. This gives us three choices for $e_2,e_3,e_4$. We then again have three choices for $e_5,e_6,e_7$, as the parity of $e_3$ is determined. Continuing this pattern gives us a total of $2(3^{n-1})$ valid MV assignments.
\end{proof}

\begin{theorem}\label{thm:2xnedges}
    The number of edges of $\OFG(M_{2,n})$ is $8(n+1)3^{n-3}$.
\end{theorem}
\begin{proof}
    We use a probabilistic method similar to one in \cite{natashaOFG}. Let $E=E(\OFG(M_{2,n}))$ be the set of edges in the origami flip graph of the $2\times n$ Miura-ori, and let $V=V(\OFG(M_{2,n}))$ be the set of vertices. Let $\mathbb{E}(\deg v)$ denote the expected degree of a randomly chosen vertex in $V$. Then
\begin{equation}\label{eq:handshake}
      \mathbb{E}(\deg v)=
\frac{\sum_{v\in V} \deg(v)}{|V|}=
\frac{2|E|}{|V|}.
    \end{equation}
Notice that $\mathbb{E}(\deg (v))=\mathbb{E}(X)$, where $X$ is the number of flippable faces of a specific MV assignment. Let $I_{i,j}$ be the indicator variable for $\alpha_{i,j}$ being flippable, so $X=I_{1,1}+...+I_{2,n}$. By linearity of expectation, 
    $$\E(X)=\sum_{i=1}^2\sum_{j=1}^{n}\E(I_{i,j}).$$
    We know that $\E(I_{i,j})=P(\alpha_{i,j} \text{ is flippable})$, and we calculate this probability in each of two cases: Any face $\alpha_{i,j}$ of $M_{2,n}$ is either bordered by three or two creases; we call these faces interior and corner respectively. 
    
    In the case where $\alpha_{i,j}$ is an interior face, its flippability is determined by the MV assignment of creases around the two vertices it is adjacent to. Without loss of generality, consider $\alpha_{1,2}$ (see Figure~\ref{fig:edgelabel}): its flippability is dependent on the MV assignment of creases around two vertices $x_1,x_2$. Thus, we consider the creases $e_2,e_3,e_5$ to determine its flippability. First, $e_2$ and $e_3$ need to be of opposite parity, since they meet at an acute angle. This occurs with probability $\frac{2}{3}$. Also, $e_3$ and $e_5$ must have the same MV parity since they meet at an obtuse angle. This also occurs with probability $\frac{2}{3}$. Since these probabilities are independent, the overall probability of flippability is $\frac{4}{9}$ for any interior face. 
    
    In the case that $\alpha_{i,j}$ is a corner face ($\alpha_{1,1},\alpha_{2,1},\alpha_{1,n},$ or $\alpha_{2,n}$), we only have one of the constraints listed above since it is only adjacent to one vertex, and thus the probability of flippability is $\frac{2}{3}$ instead of $\frac{4}{9}$. 
    
    Since we have $4$ corner faces and $2n-4$ interior faces on the $2\times n$ Miura-ori, we get 
    $$\mathbb{E}(X)=4\cdot\frac{2}{3}+(2n-4)\cdot\frac{4}{9} = \frac{8(n+1)}{9}.$$
    Plugging this value into \eqref{eq:handshake} and using the expression for $|V|$ from Theorem \ref{thm:2xnvertices} gives us 
    $$
    |E| = |V|\cdot\mathbb{E}(X)/2 = 2(3^{n-1})\E(X)/2
    = \frac{3^{n-1}(8(n+1))}{9} = 8(n+1)3^{n-3}.$$  
\end{proof}

These enumerations of vertices and edges in $\OFG(M_{2,n})$ are difficult to generalize to $M_{m,n}$. In \cite{Ginepro} there is discussion of a transfer matrix method for counting locally valid MV assignments of $M_{m,n}$ (and thus the vertices of $\OFG(M_{m,n})$), but no closed formula is known. Finding the number of edges in $\OFG(M_{m,n})$ remains an open problem.


\section{The extension forest \texorpdfstring{$\chi$}{X}}\label{subsec:extensions}

It is helpful to think about what happens when we \textit{extend} a $2 \times n$ Miura-ori to a $2 \times (n+1)$ Miura-ori by adding three creases to the right end of $M_{2,n}$ around a new vertex $x_n$. We may then refer to \textit{extending a valid MV assignment} $\mu$ of $M_{2,n}$ to a new, valid MV assignment $\mu'$ of $M_{2,n+1}$ by choosing Ms and Vs for the three new creases around $x_n$ that satisfy Maekawa's Theorem and the Big-Little-Big Lemma. This makes it significantly easier to prove several theorems about the structure of $\OFG(M_{2,n})$ with induction. To do this, we motivate the representation of these MV extensions in a \textit{directed forest}, a graph composed of several disjoint, rooted trees with directed edges. 

We encode an MV assignment of $M_{2,n}$ as a $(3n-2)$-tuple of $M$s and $V$s by taking the leftmost crease $e_0$ and recording its MV-parity in the first position, then each vertex $x_i$ and recording the MV-parity of the top, right, and bottom creases ($e_{3i-1},e_{3i},e_{3i+1}$). For example, for $M_{2,3}$, we might have  $\mu = (M,M,V,M,V,M,V)$ to denote the MV-parities of $e_0,e_2,e_3,e_4,e_5,e_6,e_7$ (though we will often omit the parentheses and commas). Let $\mu(x_i)$ be the MV-parities of the creases above, to the right of, and below the vertex $x_i$. Formally, we let $\mu(x_i)=(\mu(e_{3i-1}),\mu(e_{3i}), \mu(e_{3i+1}))$. Also, for the $M_{2,1}$ case we write $\mu(x_0)=(\mu(e_0))$. 
Extending this notation to the $m\times n$ Miura-ori, we use $\mu(x_{i,j})$ to denote the parities of the three toes of the vertex $x_{i,j}$, starting with the top, middle, then bottom crease (for $1\leq i< m$ and $1\leq j < n$).

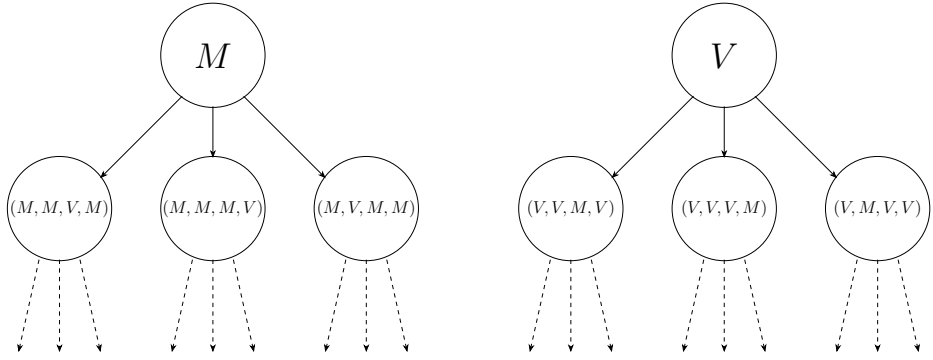
\begin{figure}[!ht]
\centering
\resizebox{0.75\textwidth}{!}{%
\begin{circuitikz}
\tikzstyle{every node}=[font=\large]
\node [font=\Huge] at (6.25,9.5) {$M$};
\node [font=\Huge] at (18.75,9.5) {$V$};
\node [font=\normalsize] at (2.5,5.75) {$(M,M,V,M)$};
\node [font=\normalsize] at (6.25,5.75) {$(M,M,M,V)$};
\node [font=\normalsize] at (10,5.75) {$(M,V,M,M)$};
\node [font=\normalsize] at (15,5.75) {$(V,V,M,V)$};
\node [font=\normalsize] at (18.75,5.75) {$(V,V,V,M)$};
\node [font=\normalsize] at (22.5,5.75) {$(V,M,V,V)$};
\draw  (2.5,5.75) circle (1.275cm);
\draw  (6.25,5.75) circle (1.275cm);
\draw  (10,5.75) circle (1.275cm);
\draw  (6.25,9.5) circle (1.275cm);
\draw  (18.75,9.5) circle (1.275cm);
\draw  (18.75,5.75) circle (1.275cm);
\draw  (22.5,5.75) circle (1.275cm);
\draw  (15,5.75) circle (1.275cm);
\draw [->, >=Stealth] (6.25,8.25) -- (6.25,7);
\draw [->, >=Stealth] (18.75,8.25) -- (18.75,7);
\draw [->, >=Stealth] (18,8.5) -- (16,6.5);
\draw [->, >=Stealth] (19.5,8.5) -- (21.5,6.5);
\draw [->, >=Stealth] (5.5,8.5) -- (3.5,6.5);
\draw [->, >=Stealth] (7,8.5) -- (9,6.5);
\draw [->, >=Stealth, dashed] (2,4.5) -- (1.5,2.25);
\draw [->, >=Stealth, dashed] (2.5,4.5) -- (2.5,2.25);
\draw [->, >=Stealth, dashed] (3,4.5) -- (3.5,2.25);
\draw [->, >=Stealth, dashed] (5.75,4.5) -- (5.25,2.25);
\draw [->, >=Stealth, dashed] (6.25,4.5) -- (6.25,2.25);
\draw [->, >=Stealth, dashed] (6.75,4.5) -- (7.25,2.25);
\draw [->, >=Stealth, dashed] (9.5,4.5) -- (9,2.25);
\draw [->, >=Stealth, dashed] (10,4.5) -- (10,2.25);
\draw [->, >=Stealth, dashed] (10.5,4.5) -- (11,2.25);
\draw [->, >=Stealth, dashed] (14.5,4.5) -- (14,2.25);
\draw [->, >=Stealth, dashed] (15,4.5) -- (15,2.25);
\draw [->, >=Stealth, dashed] (15.5,4.5) -- (16,2.25);
\draw [->, >=Stealth, dashed] (18.25,4.5) -- (17.75,2.25);
\draw [->, >=Stealth, dashed] (18.75,4.5) -- (18.75,2.25);
\draw [->, >=Stealth, dashed] (19.25,4.5) -- (19.75,2.25);
\draw [->, >=Stealth, dashed] (22,4.5) -- (21.5,2.25);
\draw [->, >=Stealth, dashed] (22.5,4.5) -- (22.5,2.25);
\draw [->, >=Stealth, dashed] (23,4.5) -- (23.5,2.25);
\end{circuitikz}
}%
\caption{The first two generations of the extension forest $\chi$, which is comprised of two infinite, directed, perfect, ternary trees (ternary meaning each vertex has $3$ children).}
\label{fig:extensionforest}
\end{figure}

Let $V_n$ be the vertex set of $\OFG(M_{2,n}).$ We define the \textit{extension forest} $\chi$ as a forest that uses the MV assignments of $V_1$ as its roots (which is either $M$ or $V$, since there is only one crease), with the second generation as valid MV assignments in $V_2$ that can be obtained via extending the original MV assignments of $V_1$, and the third generation as valid MV assignments in $V_3$ obtained via extending $V_2$, and so on. Since there are only two vertices in $V_1$, this forest is made up of two trees. We refer to the ``rows" of the forest as \textit{generations}, where the MV assignments in the $n$th generation make up $V_n$ (locally valid MV assignments of $M_{2,n}$). 
For any vertex representing an  MV assignment in $V_n$, its children are its valid extensions in $V_{n+1}$. The first $2$ generations of $\chi$ are shown in Figure \ref{fig:extensionforest}.

\begin{remark}\label{rem:1}
Each vertex of $\chi$ has $1$ parent (if it is not a root) and $3$ children.
\end{remark}

Given a vertex $u$ in the $n$th generation of $\chi$, let $f(u)$ be the number of flippable faces on the corresponding MV assignment $\mu$ of $M_{2,n}$. Note that $f(u)$ is the same as the degree of $\mu$ in $\OFG(M_{2,n})$. The following theorem serves as the motivation for the section that follows it.

\begin{theorem}\label{thm:012}
Let $u$ be a vertex in the $n$th generation of $\chi$ with corresponding MV assignment $\mu$. Let $w_1$, $w_2$, and $w_3$ be the children of $u$, with $f(w_1) \ge f(w_2) \ge f(w_3)$. If $\mu(x_{n-1})=MVM$ or $VMV$, or simply $M$ or $V$ if $n=1$, then $$(f(w_1)-f(u), f(w_2)-f(u), f(w_3)-f(u)) = (2,0,0).$$ Otherwise, $$(f(w_1)-f(u), f(w_2)-f(u), f(w_3)-f(u))=(2,1,0).$$
\end{theorem}

\begin{proof}
The result is entirely derived from the flippable face possibilities shown in Figure~\ref{flippables} when extending from $M_{2,n}$ to $M_{2,n+1}$. E.g., when
$n=1$ then there is only one crease and thus two valid MV assignments of $M_{2,1}$ ($M$ or $V$) each with two flippable faces, the extensions from which are shown in Figure \ref{fig:extensionforest}. The number of flippable faces for each MV assignment in the second generation is $4,2,2$, corresponding to adding $2,0,0$ to $f(u)$. For $n\geq 2$, if $\mu(x_{n-1})=MVM$ then there are three choices for $\mu(x_{n})$: $VMV$, $VVM$, or $MVV$. The possibilities shown in Figure~\ref{flippables} imply that if we choose $VMV$ for $\mu(x_{n})$ then we add two flippable faces to the MV assignment of $M_{2,n}$. If we choose $VMM$ or $MVV$ then we remove a flippable face from the MV assignment of $M_{2,n}$ but then add another one when extending it to $M_{2,n+1}$, for a net increase of zero flippable faces. 
Thus, if $\mu(x_{n-1})=MVM$, then $(f(w_1)-f(u), f(w_2)-f(u), f(w_3)-f(u)) = (2,0,0).$ The other possibilities for $\mu(x_{n-1})$ follow similarly from Figure~\ref{flippables}.
\end{proof}

\begin{remark}\label{rem:2}
    As can also be seen from Figure~\ref{flippables}, increasing the number of flippable faces of an MV assignment by two upon extension corresponds to assigning $\mu(x_n)$ to be $MVM$ or $VMV$.
\end{remark}


\subsection{The degree extension forest \texorpdfstring{$\chi_D$}{XD}}\label{subsec:degextension}

In addition to constructing the extension forest where vertices are labeled with MV assignments themselves, we change the labels to the number of flippable faces of the corresponding MV assignments; we define this as the \textit{degree extension forest} $\chi_D$, where the $n$th generation contains the entire degree sequence of $M_{2,n}$.  Using Theorem \ref{thm:012}, we construct $\chi_D$ starting from the first generation. To do this, we first construct the roots of each tree as two $2$s, corresponding to the MV assignments of a single mountain and a single valley. By Theorem \ref{thm:012}, we know that these two vertices have labels $4,2,2$ as their children since we add $2,0,0$. Furthermore, adding $2$ represents assigning $\mu(x_n)$ to be $MVM$ or $VMV$ (see Remark~\ref{rem:2}), so when we add $2$ to the label of a vertex $u$ to obtain one of its children $u^*$, the next set of children of $u^*$ have their labels increased by $2,0,0$. Otherwise, we add $2,1,0$ to the labels. Thus, in the third generation, the children of a vertex labeled $2$ are labeled  $4,3,2$, respectively, by adding $2,1,0$, and the children of a vertex labeled $4$ are labeled $6,4,4$ by adding $2,0,0$. 

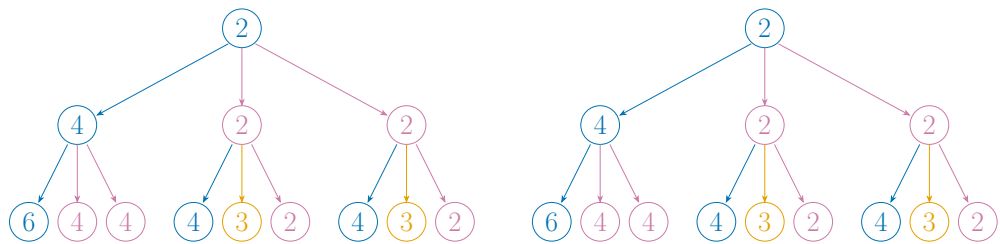
\begin{figure}
\centering
\resizebox{.8\textwidth}{!}{%
\begin{circuitikz}
\tikzstyle{every node}=[font=\LARGE]
\node [font=\LARGE, color={rgb,255:red,0; green,114; blue,178}] at (1.25,-0.5) {6};
\node [font=\LARGE, color={rgb,255:red,204; green,121; blue,167}] at (2.5,-0.5) {4};
\node [font=\LARGE, color={rgb,255:red,204; green,121; blue,167}] at (3.75,-0.5) {4};
\node [font=\LARGE, color={rgb,255:red,0; green,114; blue,178}] at (5.5,-0.5) {4};
\node [font=\LARGE, color={rgb,255:red,230; green,159; blue,0}] at (6.75,-0.5) {3};
\node [font=\LARGE, color={rgb,255:red,204; green,121; blue,167}] at (8,-0.5) {2};
\node [font=\LARGE, color={rgb,255:red,0; green,114; blue,178}] at (9.75,-0.5) {4};
\node [font=\LARGE, color={rgb,255:red,230; green,159; blue,0}] at (11,-0.5) {3};
\node [font=\LARGE, color={rgb,255:red,204; green,121; blue,167}] at (12.25,-0.5) {2};
\draw [ color={rgb,255:red,0; green,114; blue,178} ] (1.25,-0.5) circle (0.5cm);
\draw [ color={rgb,255:red,204; green,121; blue,167} ] (2.5,-0.5) circle (0.5cm);
\draw [ color={rgb,255:red,204; green,121; blue,167} ] (3.75,-0.5) circle (0.5cm);
\draw [ color={rgb,255:red,0; green,114; blue,178} ] (5.5,-0.5) circle (0.5cm);
\draw [ color={rgb,255:red,230; green,159; blue,0} ] (6.75,-0.5) circle (0.5cm);
\draw [ color={rgb,255:red,204; green,121; blue,167} ] (8,-0.5) circle (0.5cm);
\draw [ color={rgb,255:red,0; green,114; blue,178} ] (9.75,-0.5) circle (0.5cm);
\draw [ color={rgb,255:red,230; green,159; blue,0} ] (11,-0.5) circle (0.5cm);
\draw [ color={rgb,255:red,204; green,121; blue,167} ] (12.25,-0.5) circle (0.5cm);
\node [font=\LARGE, color={rgb,255:red,0; green,114; blue,178}] at (2.5,2) {4};
\node [font=\LARGE, color={rgb,255:red,204; green,121; blue,167}] at (6.75,2) {2};
\node [font=\LARGE, color={rgb,255:red,204; green,121; blue,167}] at (11,2) {2};
\node [font=\LARGE, color={rgb,255:red,0; green,114; blue,178}] at (6.75,4.5) {2};
\draw [ color={rgb,255:red,0; green,114; blue,178} ] (2.5,2) circle (0.5cm);
\draw [ color={rgb,255:red,204; green,121; blue,167} ] (6.75,2) circle (0.5cm);
\draw [ color={rgb,255:red,204; green,121; blue,167} ] (11,2) circle (0.5cm);
\draw [ color={rgb,255:red,0; green,114; blue,178} ] (6.75,4.5) circle (0.5cm);
\draw [ color={rgb,255:red,0; green,114; blue,178}, ->, >=Stealth] (6.4,4.1) -- (3,2.25);
\draw [ color={rgb,255:red,0; green,114; blue,178}, ->, >=Stealth] (2.25,1.5) -- (1.5,0);
\draw [ color={rgb,255:red,0; green,114; blue,178}, ->, >=Stealth] (6.5,1.5) -- (5.75,0);
\draw [ color={rgb,255:red,0; green,114; blue,178}, ->, >=Stealth] (10.75,1.5) -- (10,0);
\draw [ color={rgb,255:red,204; green,121; blue,167}, ->, >=Stealth] (6.75,4) -- (6.75,2.5);
\draw [ color={rgb,255:red,204; green,121; blue,167}, ->, >=Stealth] (7.1,4.1) -- (10.5,2.25);
\draw [ color={rgb,255:red,204; green,121; blue,167}, ->, >=Stealth] (7,1.5) -- (7.75,0);
\draw [ color={rgb,255:red,204; green,121; blue,167}, ->, >=Stealth] (11.25,1.5) -- (12,0);
\draw [ color={rgb,255:red,204; green,121; blue,167}, ->, >=Stealth] (2.5,1.5) -- (2.5,0);
\draw [ color={rgb,255:red,204; green,121; blue,167}, ->, >=Stealth] (2.75,1.5) -- (3.5,0);
\draw [ color={rgb,255:red,230; green,159; blue,0}, ->, >=Stealth] (6.75,1.5) -- (6.75,0);
\draw [ color={rgb,255:red,230; green,159; blue,0}, ->, >=Stealth] (11,1.5) -- (11,0);
\node [font=\LARGE, color={rgb,255:red,0; green,114; blue,178}] at (14.75,-0.5) {6};
\node [font=\LARGE, color={rgb,255:red,204; green,121; blue,167}] at (16,-0.5) {4};
\node [font=\LARGE, color={rgb,255:red,204; green,121; blue,167}] at (17.25,-0.5) {4};
\node [font=\LARGE, color={rgb,255:red,0; green,114; blue,178}] at (19,-0.5) {4};
\node [font=\LARGE, color={rgb,255:red,230; green,159; blue,0}] at (20.25,-0.5) {3};
\node [font=\LARGE, color={rgb,255:red,204; green,121; blue,167}] at (21.5,-0.5) {2};
\node [font=\LARGE, color={rgb,255:red,0; green,114; blue,178}] at (23.25,-0.5) {4};
\node [font=\LARGE, color={rgb,255:red,230; green,159; blue,0}] at (24.5,-0.5) {3};
\node [font=\LARGE, color={rgb,255:red,204; green,121; blue,167}] at (25.75,-0.5) {2};
\draw [ color={rgb,255:red,0; green,114; blue,178} ] (14.75,-0.5) circle (0.5cm);
\draw [ color={rgb,255:red,204; green,121; blue,167} ] (16,-0.5) circle (0.5cm);
\draw [ color={rgb,255:red,204; green,121; blue,167} ] (17.25,-0.5) circle (0.5cm);
\draw [ color={rgb,255:red,0; green,114; blue,178} ] (19,-0.5) circle (0.5cm);
\draw [ color={rgb,255:red,230; green,159; blue,0} ] (20.25,-0.5) circle (0.5cm);
\draw [ color={rgb,255:red,204; green,121; blue,167} ] (21.5,-0.5) circle (0.5cm);
\draw [ color={rgb,255:red,0; green,114; blue,178} ] (23.25,-0.5) circle (0.5cm);
\draw [ color={rgb,255:red,230; green,159; blue,0} ] (24.5,-0.5) circle (0.5cm);
\draw [ color={rgb,255:red,204; green,121; blue,167} ] (25.75,-0.5) circle (0.5cm);
\node [font=\LARGE, color={rgb,255:red,0; green,114; blue,178}] at (16,2) {4};
\node [font=\LARGE, color={rgb,255:red,204; green,121; blue,167}] at (20.25,2) {2};
\node [font=\LARGE, color={rgb,255:red,204; green,121; blue,167}] at (24.5,2) {2};
\node [font=\LARGE, color={rgb,255:red,0; green,114; blue,178}] at (20.25,4.5) {2};
\draw [ color={rgb,255:red,0; green,114; blue,178} ] (16,2) circle (0.5cm);
\draw [ color={rgb,255:red,204; green,121; blue,167} ] (20.25,2) circle (0.5cm);
\draw [ color={rgb,255:red,204; green,121; blue,167} ] (24.5,2) circle (0.5cm);
\draw [ color={rgb,255:red,0; green,114; blue,178} ] (20.25,4.5) circle (0.5cm);
\draw [ color={rgb,255:red,0; green,114; blue,178}, ->, >=Stealth] (19.9,4.1) -- (16.5,2.25);
\draw [ color={rgb,255:red,0; green,114; blue,178}, ->, >=Stealth] (15.75,1.5) -- (15,0);
\draw [ color={rgb,255:red,0; green,114; blue,178}, ->, >=Stealth] (20,1.5) -- (19.25,0);
\draw [ color={rgb,255:red,0; green,114; blue,178}, ->, >=Stealth] (24.25,1.5) -- (23.5,0);
\draw [ color={rgb,255:red,204; green,121; blue,167}, ->, >=Stealth] (20.25,4) -- (20.25,2.5);
\draw [ color={rgb,255:red,204; green,121; blue,167}, ->, >=Stealth] (20.6,4.1) -- (24,2.25);
\draw [ color={rgb,255:red,204; green,121; blue,167}, ->, >=Stealth] (20.5,1.5) -- (21.25,0);
\draw [ color={rgb,255:red,204; green,121; blue,167}, ->, >=Stealth] (24.75,1.5) -- (25.5,0);
\draw [ color={rgb,255:red,204; green,121; blue,167}, ->, >=Stealth] (16,1.5) -- (16,0);
\draw [ color={rgb,255:red,204; green,121; blue,167}, ->, >=Stealth] (16.25,1.5) -- (17,0);
\draw [ color={rgb,255:red,230; green,159; blue,0}, ->, >=Stealth] (20.25,1.5) -- (20.25,0);
\draw [ color={rgb,255:red,230; green,159; blue,0}, ->, >=Stealth] (24.5,1.5) -- (24.5,0);
\end{circuitikz}
}%
\caption{The first $3$ generations of the degree extension forest $\chi_D$.}
\label{fig:degextensionforest}
\end{figure}

Tracking whether or not we add $2,0,0$ or $2,1,0$ to a vertex when constructing $\chi_D$ is not any easier than constructing the extension forest $\chi$, since we have to track what MV assignment corresponds to each vertex. To better track whether or not a degree extends by adding $2,1,0$ or $2,0,0$ (by Theorem \ref{thm:012}), we color the edges between generations and the vertices themselves either magenta, orange, or blue. Edges where the difference between the child and parent is $0$ are colored magenta; edges where the difference is $1$ are colored orange; edges where the difference is $2$ are colored blue. The children are colored the same as the edge that joins them with their parent. By Remark~\ref{rem:2}, if a vertex is blue, meaning we added $2$ flippable faces, then its corresponding MV assignment ends in $MVM$ or $VMV$.  Then, by Theorem~\ref{thm:012} itself, we add $2,0,0$ to its label obtain its children. Otherwise, if a vertex is orange or magenta, we add $2,1,0$ to it. This coloring approach allows us to construct $\chi_D$ without tracking the actual MV assignments that these vertices correspond to, only the colors. The first three generations of $\chi_D$ are shown in Figure \ref{fig:degextensionforest}.\\

\begin{remark}\label{rem:3}
    Every parent vertex $u \in \chi_D$ has at least one child joined by a magenta edge (adding $0$), and exactly one child joined by a blue edge (adding $2$).
\end{remark}

We will prove many facts about the degree sequence of $\OFG(M_{2,n})$ using the degree extension forest $\chi_D$, both in this section and in Section \ref{subsec:polynomial}.

\begin{lemma}\label{thm:oppositeexist}
    In the $n$th generation of $\chi_D$, the minimum labeling is $2$, and for all natural $n \ge 2$, it appears four times and is colored magenta.
\end{lemma}

\begin{proof}
    In the first generation of $\chi_D$, there are two blue vertices, each labeled with $2$. Labelings cannot decrease from generation to generation, but they can remain constant by adding $0$. Thus, the minimum labeling is $2$ in any generation. In the second generation, there are four magenta $2$s. Since we can only add $0$ to the parent exactly once by Theorem \ref{thm:012}, there will be exactly four $2$s in all subsequent generations. Adding $0$ also means that each of these $2$s will be magenta, thus completing the proof.
\end{proof}

\begin{figure}
\centering
\resizebox{0.65\textwidth}{!}{%
\begin{circuitikz}
\tikzstyle{every node}=[font=\Huge]
\draw [short] (7.5,7) -- (7,5.75);
\draw [short] (7.5,7) -- (10,7);
\draw [short] (10,7) -- (9.5,5.75);
\draw [short] (9.5,5.75) -- (10,4.5);
\draw [short] (7,5.75) -- (7.5,4.5);
\draw [short] (7.5,4.5) -- (10,4.5);
\draw [dashed] (10.25,5.75) -- (11.25,5.75);
\draw [short] (12.25,7) -- (11.75,5.75);
\draw [short] (12.25,7) -- (14.75,7);
\draw [short] (14.75,7) -- (14.25,5.75);
\draw [short] (14.25,5.75) -- (14.75,4.5);
\draw [short] (11.75,5.75) -- (12.25,4.5);
\draw [short] (12.25,4.5) -- (14.75,4.5);
\draw [short] (7.5,3.75) -- (7,2.5);
\draw [short] (7.5,3.75) -- (10,3.75);
\draw [short] (10,3.75) -- (9.5,2.5);
\draw [short] (9.5,2.5) -- (10,1.25);
\draw [short] (7,2.5) -- (7.5,1.25);
\draw [short] (7.5,1.25) -- (10,1.25);
\draw [dashed] (10.25,2.5) -- (11.25,2.5);
\draw [short] (12.25,3.75) -- (11.75,2.5);
\draw [short] (12.25,3.75) -- (14.75,3.75);
\draw [short] (14.75,3.75) -- (14.25,2.5);
\draw [short] (14.25,2.5) -- (14.75,1.25);
\draw [short] (11.75,2.5) -- (12.25,1.25);
\draw [short] (12.25,1.25) -- (14.75,1.25);
\draw [short] (17.25,7) -- (16.75,5.75);
\draw [short] (17.25,7) -- (19.75,7);
\draw [short] (19.75,7) -- (19.25,5.75);
\draw [short] (19.25,5.75) -- (19.75,4.5);
\draw [short] (16.75,5.75) -- (17.25,4.5);
\draw [short] (17.25,4.5) -- (19.75,4.5);
\draw [dashed] (20,5.75) -- (21,5.75);
\draw [short] (22,7) -- (21.5,5.75);
\draw [short] (22,7) -- (24.5,7);
\draw [short] (24.5,7) -- (24,5.75);
\draw [short] (24,5.75) -- (24.5,4.5);
\draw [short] (21.5,5.75) -- (22,4.5);
\draw [short] (22,4.5) -- (24.5,4.5);
\draw [short] (17.25,3.75) -- (16.75,2.5);
\draw [short] (17.25,3.75) -- (19.75,3.75);
\draw [short] (19.75,3.75) -- (19.25,2.5);
\draw [short] (19.25,2.5) -- (19.75,1.25);
\draw [short] (16.75,2.5) -- (17.25,1.25);
\draw [short] (17.25,1.25) -- (19.75,1.25);
\draw [dashed] (20,2.5) -- (21,2.5);
\draw [short] (22,3.75) -- (21.5,2.5);
\draw [short] (22,3.75) -- (24.5,3.75);
\draw [short] (24.5,3.75) -- (24,2.5);
\draw [short] (24,2.5) -- (24.5,1.25);
\draw [short] (21.5,2.5) -- (22,1.25);
\draw [short] (22,1.25) -- (24.5,1.25);
\draw [line width=1pt, short] (7,5.75) -- (9.5,5.75);
\draw [line width=1pt, short] (8.25,5.75) -- (8.75,7);
\draw [line width=1pt, short] (11.75,5.75) -- (14.25,5.75);
\draw [line width=1pt, short] (13.5,7) -- (13,5.75);
\draw [line width=1pt, short] (16.75,5.75) -- (19.25,5.75);
\draw [line width=1pt, short] (21.5,5.75) -- (24,5.75);
\draw [line width=1pt, short] (18,5.75) -- (18.5,4.5);
\draw [line width=1pt, short] (22.75,5.75) -- (23.25,4.5);
\draw [line width=1pt, short] (8.25,2.5) -- (8.75,1.25);
\draw [line width=1pt, short] (13,2.5) -- (13.5,1.25);
\draw [line width=1pt, short] (18.5,3.75) -- (18,2.5);
\draw [line width=1pt, short] (23.25,3.75) -- (22.75,2.5);
\draw [line width=1pt, dashed] (8.25,5.75) -- (8.75,4.5);
\draw [line width=1pt, dashed] (13,5.75) -- (13.5,4.5);
\draw [line width=1pt, dashed] (7,2.5) -- (9.5,2.5);
\draw [line width=1pt, dashed] (8.25,2.5) -- (8.75,3.75);
\draw [line width=1pt, dashed] (11.75,2.5) -- (14.25,2.5);
\draw [line width=1pt, dashed] (13,2.5) -- (13.5,3.75);
\draw [line width=1pt, dashed] (16.75,2.5) -- (19.25,2.5);
\draw [line width=1pt, dashed] (18,2.5) -- (18.5,1.25);
\draw [line width=1pt, dashed] (21.5,2.5) -- (24,2.5);
\draw [line width=1pt, dashed] (22.75,2.5) -- (23.25,1.25);
\draw [line width=1pt, dashed] (22.75,5.75) -- (23.25,7);
\draw [line width=1pt, dashed] (18,5.75) -- (18.5,7);
\node [font=\Huge] at (6.25,5.75) {$\nu_1$};
\node [font=\Huge] at (6.25,2.5) {$\nu_2$};
\node [font=\Huge] at (16,5.75) {$\nu_3$};
\node [font=\Huge] at (16,2.5) {$\nu_4$};
\end{circuitikz}
}%
\caption{The explicit constructions of the $4$ degree-$2$ vertices of $\OFG(M_{2,n})$. For $\nu_1$ and $\nu_2$, only $\alpha_{1,1}$ and $\alpha_{2,n}$ are flippable, and for $\nu_3$ and $\nu_4$, only $\alpha_{2,1}$ and $\alpha_{1,n}$ are flippable. (Dots indicate the MV assignment repeating $n-1$ times.)}
\label{fig:deg2construction}
\end{figure}
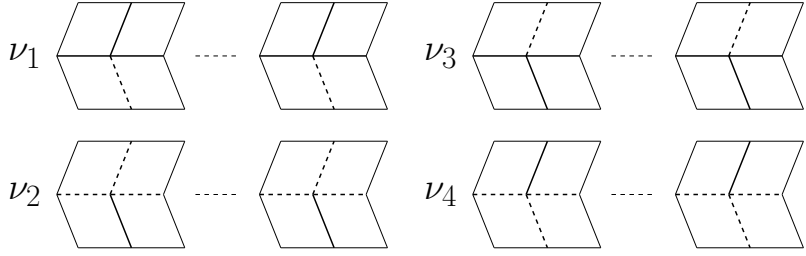

Lemma \ref{thm:oppositeexist} is interesting since it implies that the number of degree-$2$ vertices in $\OFG(M_{2,n})$ is constant as $n$ increases despite the number of vertices increasing exponentially. The four MV assignments corresponding to these degree-$2$ vertices are constructed in Figure \ref{fig:deg2construction}. 
In fact, we find that this result holds for $\OFG(M_{m,n})$ as well.
Let us denote the crease between faces $\alpha_{i,j}$ and $\alpha_{i',j'}$ by $\{\alpha_{i,j}, \alpha_{i',j'}\}$ (here either $i'=i$ and $j'=j+1$ or $i'=i+1$ and $j'=j$).

\begin{theorem}\label{thm:4verticesdeg2}
    In $\OFG(M_{m,n})$, there exist exactly $4$ vertices of degree 2. 
\end{theorem}
\begin{proof}
    We construct four MV assignments that only have two flippable faces, and then show that these are the only such MV assignments possible.
Consider the MV assignment on $M_{m,n}$ where creases $\{\alpha_{i,j}, \alpha_{i,j+1}\}$ for $1\leq i\leq m-1$ odd and any $1\leq j\leq n-1$ are mountain and all other creases are valley (the $M_{3,4}$ example of this is shown in the upper-left crease pattern in Figure~\ref{fig:deg2vertices}).  This MV assignment has only two flippable faces: $\alpha_{1,n}$ and $\alpha_{m,1}$. To see this is the case, consider any $\alpha_{i,j}$. Flipping $\alpha_{i,j}$ makes the vertex $x_{i,j}$ or the vertex $x_{i-1,j-1}$, if they exist, violate Theorem \ref{biglittlebig}, since the three toes of the vertex will all be the same. $x_{i-1,j-1}$ or $x_{i,j}$ exist for all $(i,j) \neq (1,n), (m,1)$, so none of those $\alpha_{i,j}$ faces are flippable. For the face $\alpha_{1,n}$, the only vertex it borders is $x_{1,n-1}$ and the only creases it borders are $\{\alpha_{1,n-1}, \alpha_{1,n}\}$ and $\{\alpha_{1,n},\alpha_{2,n}\}$, which have different MV parities. Therefore if we flip $\alpha_{1,n}$ the MV assignment at $x_{1,n-1}$ will still be valid since the two creases simply swap parities. Flipping $\alpha_{m,1}$ is similar: if $m$ is odd, creases $\{\alpha_{m-1,1},\alpha_{m,1}\}$ and $\{\alpha_{m,1},\alpha_{m,2}\}$ will swap parities, and if $m$ is even then $\{\alpha_{m-1,1},\alpha_{m,1}\}$ and $\{\alpha_{m,1},\alpha_{m,2}\}$ switch from being both valley to both mountain, which turns $x_{m-1,1}$ into a majority mountain vertex. A similar analysis can be done to show that the MV assignment (shown in the upper-right of Figure~\ref{fig:deg2vertices}) where creases $\{\alpha_{i,j}, \alpha_{i,j+1}\}$ for $1\leq i\leq m-1$ even and any $1\leq j\leq n-1$ are mountain and all other creases are valley has only $\alpha_{1,1}$ and $\alpha_{m,n}$ flippable.  To get the other two MV assignments, we can simply take these two MV assignments and swap the parities of all the creases (as seen in the bottom row of Figure \ref{fig:deg2vertices}). 

\begin{figure}
\centering
\resizebox{0.8\linewidth}{!}{%
\begin{circuitikz}
\tikzstyle{every node}=[font=\LARGE]

\begin{scope}[xshift=0cm,yshift=0cm]
\draw [short] (11.25,9.75) -- (10,7.25);
\draw [dashed] (10,7.25) -- (11.25,4.75);
\draw [dashed] (10,7.25) -- (13.75,7.25);
\draw [short] (13.75,7.25) -- (15,9.5);
\draw [dashed] (13.75,7.25) -- (17.5,7.25);
\draw [dashed] (13.75,7.25) -- (15,4.75);
\draw [short] (18.75,9.75) -- (17.5,7.25);
\draw [dashed] (17.5,7.25) -- (21.25,7.25);
\draw [dashed] (17.5,7.25) -- (18.75,4.75);
\draw [dashed] (18.75,4.75) -- (21.25,4.75);
\draw [short] (18.75,4.75) -- (17.5,2.25);
\draw [dashed] (15,4.75) -- (18.75,4.75);
\draw [dashed] (11.25,4.75) -- (15,4.75);
\draw [short] (15,4.75) -- (13.75,2.25);
\draw [dashed] (7.5,7.25) -- (10,7.25);
\draw [dashed] (7.5,4.75) -- (11.25,4.75);
\draw [short] (11.25,4.75) -- (10,2.25);
\node [font=\Huge, color={rgb,255:red,0; green,200; blue,50}] at (20,8.5) {$\checkmark$};
\node [font=\Huge, color={rgb,255:red,0; green,200; blue,50}] at (8.75,3.5) {$\checkmark$};
\end{scope}

\begin{scope}[xshift=16cm,yshift=0cm]
\begin{scope}
\draw [dashed] (11.25,9.75) -- (10,7.25);
\draw [short] (10,7.25) -- (11.25,4.75);
\draw [dashed] (10,7.25) -- (13.75,7.25);
\draw [dashed] (13.75,7.25) -- (15,9.5);
\draw [dashed] (13.75,7.25) -- (17.5,7.25);
\draw [short] (13.75,7.25) -- (15,4.75);
\draw [dashed] (18.75,9.75) -- (17.5,7.25);
\draw [dashed] (17.5,7.25) -- (21.25,7.25);
\draw [short] (17.5,7.25) -- (18.75,4.75);
\draw [dashed] (18.75,4.75) -- (21.25,4.75);
\draw [dashed] (18.75,4.75) -- (17.5,2.25);
\draw [dashed] (15,4.75) -- (18.75,4.75);
\draw [dashed] (11.25,4.75) -- (15,4.75);
\draw [dashed] (15,4.75) -- (13.75,2.25);
\draw [dashed] (7.5,7.25) -- (10,7.25);
\draw [dashed] (7.5,4.75) -- (11.25,4.75);
\draw [dashed] (11.25,4.75) -- (10,2.25);
\node [font=\Huge, color={rgb,255:red,0; green,200; blue,50}] at (20,3.5) {$\checkmark$};
\node [font=\Huge, color={rgb,255:red,0; green,200; blue,50}] at (8.75,8.5) {$\checkmark$};
\end{scope}
\end{scope}

\begin{scope}[xshift=0cm,yshift=-9cm]
\begin{scope}
\draw [dashed] (11.25,9.75) -- (10,7.25);
\draw [short] (10,7.25) -- (11.25,4.75);
\draw [short] (10,7.25) -- (13.75,7.25);
\draw [dashed] (13.75,7.25) -- (15,9.5);
\draw [short] (13.75,7.25) -- (17.5,7.25);
\draw [short] (13.75,7.25) -- (15,4.75);
\draw [dashed] (18.75,9.75) -- (17.5,7.25);
\draw [short] (17.5,7.25) -- (21.25,7.25);
\draw [short] (17.5,7.25) -- (18.75,4.75);
\draw [short] (18.75,4.75) -- (21.25,4.75);
\draw [dashed] (18.75,4.75) -- (17.5,2.25);
\draw [short] (15,4.75) -- (18.75,4.75);
\draw [short] (11.25,4.75) -- (15,4.75);
\draw [dashed] (15,4.75) -- (13.75,2.25);
\draw [short] (7.5,7.25) -- (10,7.25);
\draw [short] (7.5,4.75) -- (11.25,4.75);
\draw [dashed] (11.25,4.75) -- (10,2.25);
\node [font=\Huge, color={rgb,255:red,0; green,200; blue,50}] at (20,8.5) {$\checkmark$};
\node [font=\Huge, color={rgb,255:red,0; green,200; blue,50}] at (8.75,3.5) {$\checkmark$};
\end{scope}
\end{scope}

\begin{scope}[xshift=16cm,yshift=-9cm]
\begin{scope}
\draw [short] (11.25,9.75) -- (10,7.25);
\draw [dashed] (10,7.25) -- (11.25,4.75);
\draw [short] (10,7.25) -- (13.75,7.25);
\draw [short] (13.75,7.25) -- (15,9.5);
\draw [short] (13.75,7.25) -- (17.5,7.25);
\draw [dashed] (13.75,7.25) -- (15,4.75);
\draw [short] (18.75,9.75) -- (17.5,7.25);
\draw [short] (17.5,7.25) -- (21.25,7.25);
\draw [dashed] (17.5,7.25) -- (18.75,4.75);
\draw [short] (18.75,4.75) -- (21.25,4.75);
\draw [short] (18.75,4.75) -- (17.5,2.25);
\draw [short] (15,4.75) -- (18.75,4.75);
\draw [short] (11.25,4.75) -- (15,4.75);
\draw [short] (15,4.75) -- (13.75,2.25);
\draw [short] (7.5,7.25) -- (10,7.25);
\draw [short] (7.5,4.75) -- (11.25,4.75);
\draw [short] (11.25,4.75) -- (10,2.25);
\node [font=\Huge, color={rgb,255:red,0; green,200; blue,50}] at (20,3.5) {$\checkmark$};
\node [font=\Huge, color={rgb,255:red,0; green,200; blue,50}] at (8.75,8.5) {$\checkmark$};
\end{scope}
\end{scope}

\end{circuitikz}
}
\caption{The four MV assignments that have only two flippable faces.}
\label{fig:deg2vertices}
\end{figure}

    To show these are the only four MV assignments that have two flippable faces, consider the faces along the four sides of the Miura-ori, i.e. the $4$ sets $ \{\alpha_{1,1}, \alpha_{1,2}, ... , \alpha_{1,n}\}$, $\{\alpha_{1,1}, \alpha_{2,1}, ... , \alpha_{m,1}\}$, $\{\alpha_{1,n}, \alpha_{2,n}, ... , \alpha_{m,n}\}$, and $ \{\alpha_{m,1}, \alpha_{m,2}, ... , \alpha_{m,n}\}$. We prove the following claim: each side must have at least one flippable face. 
    
    Consider $\alpha_{1,1}$. If it is flippable, then $ \{\alpha_{1,1}, \alpha_{1,2}, ... , \alpha_{1,n}\}$ and $\{\alpha_{1,1}, \alpha_{2,1}, ... , \alpha_{m,1}\}$ must both have at least $1$ flippable face. Else, it is unflippable and we know the MV assignment around $x_{1,1}$ must be either $\{\alpha_{1,1},\alpha_{2,1}\} = M$, $\mu(x_{1,1}) = (V,M,M)$, or $\{\alpha_{1,1},\alpha_{2,1}\} = V$, $\mu(x_{1,1}) = (M,V,V)$. Without loss of generality, suppose the latter.  
 
 Then consider the crease  $c=\{\alpha_{1,2},\alpha_{1,3}\}$. If $c$ is a valley, then $\alpha_{1,2}$ is flippable, since flipping it makes $\{\alpha_{1,1},\alpha_{1,2}\}$ and $\{\alpha_{1,2},\alpha_{2,2}\}$ swap MV parity (and thus $x_{1,1}$ is valid) and makes $\{\alpha_{1,2},\alpha_{2,2}\}$ and $c$ switch from VV to MM (thus $x_{1,2}$ is valid). If $c$ is a mountain, then $\alpha_{1,2}$ is unflippable and we know that the MV assignment around $x_{1,2}$ is $\mu(x_{1,2})=(M,V,V)$ (see Figure \ref{fig:deg2flippable}). Repeating this, we get either a flippable face among $\alpha_{1,3},\ldots, \alpha_{1,n-1}$ or all of $\alpha_{1,1},\ldots,\alpha_{1,n-1}$ are unflippable and $\mu(x_{1,n-1})=(M,V,V)$. But then $\alpha_{1,n}$ is flippable, since doing so swaps the parity of $\{\alpha_{1,n-1},\alpha_{1,n}\}$ and $\{\alpha_{1,n},\alpha_{2,n}\}$. Thus, $ \{\alpha_{1,1}, \alpha_{1,2}, \ldots , \alpha_{1,n}\}$ has a flippable face. 

 \begin{figure}
\centering
\resizebox{0.8\linewidth}{!}{%
\begin{circuitikz}
\tikzstyle{every node}=[font=\LARGE]

\begin{scope}[xshift=0cm,yshift=0cm]
\draw [short] (11.25,9.75) -- (10,7.25);
\draw [dashed] (10,7.25) -- (11.25,4.75);
\draw [dashed] (10,7.25) -- (13.75,7.25);
\draw [dashed] (13.75,7.25) -- (15,9.5);
\draw [dashed] (7.5,7.25) -- (10,7.25);

\node [font=\LARGE, color={rgb,255:red,0; green,200; blue,50}] at (12.5,8.5) {$\alpha_{1,2} \ \checkmark$};
\end{scope}

\begin{scope}[xshift=10cm,yshift=0cm]
\begin{scope}
\draw [short] (11.25,9.75) -- (10,7.25);
\draw [dashed] (10,7.25) -- (11.25,4.75);
\draw [dashed] (10,7.25) -- (13.75,7.25);
\draw [short] (13.75,7.25) -- (15,9.5);
\draw [dashed] (13.75,7.25) -- (17.5,7.25);
\draw [dashed] (13.75,7.25) -- (15,4.75);
\draw [dashed] (7.5,7.25) -- (10,7.25);
\node [font=\LARGE, color={rgb,255:red,255; green,0; blue,0}] at (12.5,8.5) {$\alpha_{1,2} \ \times$};
\end{scope}
\end{scope}

\end{circuitikz}
}
\caption{If $x_{1,1} = (M,V,V)$, then $\alpha_{1,2}$ is only flippable if $\{\alpha_{1,2},\alpha_{1,3}\}$ is a valley crease}
\label{fig:deg2flippable}
\end{figure}
 
We now examine the column of faces $\{\alpha_{1,1}, \alpha_{2,1}, \ldots, \alpha_{m,1}\}$, assuming $\alpha_{1,1}$ is not flippable and $\mu(x_{1,1})=(M,V,V)$.  If crease $c=\{\alpha_{2,1},\alpha_{3,1}\}$ is a mountain, then $\alpha_{2,1}$ is flippable and $\{\alpha_{1,1}, \alpha_{2,1}, \ldots , \alpha_{m,1}\}$ has a flippable face. If $\alpha_{2,1}$ is not flippable, then the MV assignment around $x_{2,1}$ must be: $c$ is a valley and $\mu(x_{2,1})=(V,V,M)$.

We can repeat a similar analysis for $\alpha_{3,1},\ldots,\alpha_{m-1,1}$. If $\alpha_{2,1},...,\alpha_{m-1,1}$ are all unflippable, then we 
 have that the creases around $x_{m-1,1}$ are, if $m$ is even, $\mu(x_{m-1,1})=(M,V,V)$ and the other, left crease valley, and if $m$ is odd we have $\mu(x_{m-1,1})=(V,V,M)$ and the left crease valley. In both cases flipping $\alpha_{m,1}$ maintains validity. Thus, $\{\alpha_{1,1}, \alpha_{2,1}, \ldots , \alpha_{m,1}\}$ has a flippable face.
 
 A similar argument can be done for the faces $\{\alpha_{1,n}, \alpha_{2,n}, \ldots , \alpha_{m,n}\}$ and $ \{\alpha_{m,1}, \alpha_{m,2}, \ldots , \alpha_{m,n}\}$. 

    Thus, we have that each side must have at least one flippable face. This means that the minimum number of flippable faces is 2, where the only flippable faces are in the two opposite corners, which corresponds to the 4 MV assignments we constructed above. 
\end{proof}

We now return to the degree extension forest $\chi_D$ and the $M_{2,n}$ case.

\begin{lemma}\label{lem:2nblue}
    For all $n\in\N$, the maximum labeling in the $n$th generation of $\chi_D$ is $2n$, appears twice, and is colored blue.
\end{lemma}

\begin{proof}
    When $n=1$, the maximum $2n=2$ appears twice and is always blue. The maximum of the $(n+1)$th generation is generated by adding $2$ to the maximum of the $n$th generation, hence, the next maximum is $2(n+1)$ and is also blue. Furthermore, a parent can only add $2$ to its labeling for exactly one of its children, so the number of vertices with labels of $2n$ does not change across generations. Thus, the maximum in the $n$th generation is $2n$, appears twice, and is always blue. 
\end{proof}

\begin{remark}\label{rem:4}
    The first generation in which an even number $2a$ appears in $\chi_D$ is the $a$th generation, where $2a$ is the maximum of the labelings in that generation.
\end{remark}

\begin{lemma}\label{lem:2n-1}
    For all $n\in \N$, no vertices are labeled $2n-1$ in the $n$th generation of $\chi_D$.
\end{lemma}

\begin{proof}
    We prove this via induction on $n$. The base case $n=1$ is evident in Figure \ref{fig:degextensionforest}, since $1$ is not in the $1$st generation. For the inductive step, 
    assume that $2k-1$ is not in the $k$th generation. If we were to generate the number $2k+1$ in the $(k+1)$th generation, then in the $k$th generation we would need to either add 1 or 2 to some $2k-1$ label, which the induction hypothesis prohibits, or add $0$ to some $2k+1$, which cannot happen by  Lemma \ref{lem:2nblue}. Thus, $2k+1$ is not in the $(k+1)$th generation if $2k-1$ is not in the $k$th generation.
\end{proof}

Let $\Set(n)$ be a function that inputs a natural number $n$ and outputs the set of all numeric labelings of the vertices in the $n$th generation of $\chi_D$ (meaning duplicate entries are removed). This equivalently gives the possible degrees of vertices in $\OFG(M_{2,n})$, which we determine using the following theorem.

\begin{theorem}\label{thm:set}
    For all $n \ge 2$, $\Set(n)=\{2,3,4,...,2n-2,2n\}$.
\end{theorem}

\begin{proof}
By the previous Lemmas we know   $\min(\Set(n))=2$, $\max(\Set(n))=2n$, and $2n-1\not\in\Set(n)$. Thus we only need to show that all other naturals between 2 and $2n$ are in the set. 

Figure \ref{fig:degextensionforest} shows that $\Set(2)=\{2,4\}$ and $\Set(3)=\{2,3,4,6\}$, so the Theorem is true for $n=2$ and $n=3$. Proceeding by induction, assume for $k\geq3$ that $\Set(k)=\{2,3,4,...,2k-2,2k\}$. We know by the previous Lemmas that $2k+1 \notin \Set(k+1)$ and $\max(\Set(k+1)) = 2k+2$. Since all vertices have a child with equal numeric label (the magenta-colored child), we know that $\{2,3,4,...,2k-2,2k\} \subset \Set(k+1)$. All that is left to show is that $2k-1\in \Set(k+1)$, but this is true because $2k-3\in\Set(k)$ and each parent has a blue-colored child, which adds 2 to its label to get $2k-1$. Since 2 is our minimum and $2k+2$ is our maximum, there are no other elements in $\Set(k+1)$, and we are done.

\end{proof}


\subsection{Number of vertices of degree \texorpdfstring{$d$}{d} in \texorpdfstring{$\OFG(M_{2,n})$}{OFG(M2n)}}\label{subsec:polynomial}

In order to search for patterns in the degree sequences for $\OFG(M_{2,n})$, we wrote code to calculate the number of flippable faces for all MV assignments of $M_{2,n}$ for small $n$.
Table \ref{tab:degdistribution} shows the results, where $v_n^d$ denotes the number of vertices in $\OFG(M_{2,n})$ with degree $d$, i.e., the number of vertices in the $n$th generation of $\chi_D$ with label $d$.

\begin{table}[h]
\centering
\begin{tabular}{c|c|c|c|c|c|c|c|c}
\textbf{$d$} & \textbf{$v_2^d$} & \textbf{$v_3^d$} & \textbf{$v_4^d$} & \textbf{$v_5^d$} & \textbf{$v_6^d$} & \textbf{$v_7^d$} & \textbf{$v_8^d$} & \textbf{$v_9^d$} \\ \hline
\textbf{2} & 4 & 4 & 4 & 4 & 4 & 4 & 4 & 4 \\ \hline
\textbf{3} & 0 & 4 & 8 & 12 & 16 & 20 & 24 & 28\\ \hline
\textbf{4} & 2 & 8 & 20 & 36 & 56 & 80 & 108 & 140\\ \hline
\textbf{5} &  & 0 & 8 & 36 & 88 & 168 & 280 & 428\\ \hline
\textbf{6} &  & 2 & 12 & 44 & 128 & 296 & 584 & 1032\\ \hline
\textbf{7} &  &  & 0 & 12 & 80 & 296 & 792 & 1744\\ \hline
\textbf{8} &  &  & 2 & 16 & 76 & 292 & 924 & 2428\\ \hline
\textbf{9} &  &  &  & 0 & 16 & 140 & 680 & 2396\\ \hline
\textbf{10} &  &  &  & 2 & 20 & 116 & 544 & 2144\\ \hline
\textbf{11} &  &  &  &  & 0 & 20 & 216 & 1288\\ \hline
\textbf{12} &  &  &  &  & 2 & 24 & 164 & 900\\ \hline
\textbf{13} &  &  &  &  &  & 0 & 24 & 308\\ \hline
\textbf{14} &  &  &  &  &  & 2 & 28 & 220\\ \hline
\textbf{15} &  &  &  &  &  &  & 0 & 28\\ \hline
\textbf{16} &  &  &  &  &  &  & 2 & 32\\ \hline
\textbf{17} &  &  &  &  &  &  &  & 0\\ \hline
\textbf{18} &  &  &  &  &  &  &  & 2\\
\end{tabular}

\caption{Degree distributions of $\OFG(M_{2,n})$ for $2 \le n \le 9.$ Empty cells are $0$.}
\label{tab:degdistribution}
\end{table}

By reading each row of Table \ref{tab:degdistribution}, it appears that each forms a sequence described by polynomials of increasing degree: the $d=2$ row is constant, the $d=3$ row is linear, the $d=4$ row is quadratic, and so on. We seek to prove this via a recursive technique, but first we construct recursive equations for $v_n^d$ by using properties of the degree extension forest $\chi_D$.\\

\begin{proposition}\label{prop:recurrences}
    Let $v_n^d$ be the number of vertices in the $n$th generation of $\chi_D$ labeled $d$. Let $b_n^d$ be the number of such vertices that are blue, and $w_n^d$ be the number of such vertices that are not blue (W for warm). Then these variables satisfy  the following recursions:
    $$b_n^d=v_{n-1}^{d-2},$$ $$w_n^d=v_{n-1}^d+w_{n-1}^{d-1}+v_{n-2}^{d-2}, \text{and}$$ $$v_{n}^d=v_{n-1}^d+w_{n-1}^{d-1}+v_{n-1}^{d-2}+v_{n-2}^{d-2}.$$
\end{proposition}

\begin{proof}
    The first recursion follows from Theorem~\ref{thm:012} and Remark~\ref{rem:3}.

    We now consider the equation for $w_n^d$. A non-blue vertex in the $n$th generation labeled with $d$ must have its parent of the $(n-1)$th generation be a vertex in the $(n-1)$th generation labeled with either $d$ (any color) or $d-1$ (non-blue). However, if the parent is labeled $d$, then there are two cases: if the parent is non-blue, it generates one child of degree $d$, and if the parent is blue, it generates two children of degree $d$, by Theorem \ref{thm:012}. Thus, we get that $w_n^d=w_{n-1}^d+2b_{n-1}^d+w_{n-1}^{d-1}.$ However, since $w_{n-1}^d+b_{n-1}^d=v_{n-1}^d$ and $b_{n-1}^d=v_{n-2}^{d-2}$, we rewrite this expression as $w_n^d=v_{n-1}^d+w_{n-1}^{d-1}+v_{n-2}^{d-2}$, giving us the second equation. 
    
    We know that $v_n^d=b_n^d+w_n^d$, so adding the first and second equations gives us the third.
\end{proof}

With the recursive equations from Proposition \ref{prop:recurrences}, we prove that $v_n^d$ is a degree-$(d-2)$ polynomial in terms of $n$. We begin by proving some smaller cases to provide intuition for the proof. We have already proven that the explicit formula describing $v_n^2$ is a degree-$0$ polynomial in Lemma~\ref{thm:oppositeexist} (note that this excludes the case $n=1$ where there are $2$ vertices of degree $2$). Consider the case of $d=3$; we show that $v_n^3$ is described by a degree-$1$ polynomial. Since a vertex with the labeling $3$ can have its parent be $2$ or $3$ (always non-blue, with only one child labeled $3$), we know that $$v_n^3=v_{n-1}^2+v_{n-1}^3=4+v_{n-1}^3.$$ 
This recursion solves to the first-degree polynomial
$v_n^3=4(n-2)$.  We now generalize this technique for all natural $d \ge 2$ using a strong induction argument.

\begin{theorem}\label{thm:polynomials}
    For all natural $d \ge 2$, the explicit formula for $v_n^d$ where $n \ge \left\lceil \frac{d}{2} \right\rceil+1$ is a polynomial of $n$ with degree $d-2$.
\end{theorem}

\begin{proof}
    Setting $n \ge \left\lceil \frac{d}{2} \right\rceil+1$ includes the non-zero values of $v_n^d$ but excludes the edge case of $v_{d/2}^d=2$ for $d$ even.  We prove the statement by strong induction on $d$. We have already proven the base cases $d=2,3$, so we proceed to the inductive step. 

    Assume that for all $d \le k$, the explicit formula for $v_n^d$ is a polynomial of $n$ with degree $(d-2)$. We show this implies that $v_n^{k+1}$ is described by a degree-$(k-1)$ polynomial. By Proposition \ref{prop:recurrences}, we have $$v_n^{k+1}=v_{n-1}^{k+1}+w_{n-1}^k+v_{n-1}^{k-1}+v_{n-2}^{k-1}.$$ Again, we use Proposition \ref{prop:recurrences} to replace the $w_{n-1}^k$ term, giving us that 
    \begin{equation}\label{eq:vnk+1}
    v_n^{k+1}=v_{n-1}^{k+1}+v_{n-2}^k+w_{n-2}^{k-1}+v_{n-3}^{k-2}+v_{n-1}^{k-1}+v_{n-2}^{k-1}.
    \end{equation}
    Since $w_{n-2}^{k-1}=v_{n-2}^{k-1}-b_{n-2}^{k-1}=v_{n-2}^{k-1}-v_{n-3}^{k-3}$ (using Proposition~\ref{prop:recurrences}), the inductive hypothesis implies that $w_{n-2}^{k-1}$ is described by a polynomial of at most degree $k-3$. In Equation~\eqref{eq:vnk+1} the $v_{n-2}^k$ term is described, via our inductive hypothesis, by a degree-$(k-2)$ polynomial. Everything after this term provides lower-order terms to the polynomial and thus does not increase the degree. So, we simplify our expression to $$v_n^{k+1}=v_{n-1}^{k+1}+(\text{some degree-}(k-2)\text{ polynomial}).$$ 
    Telescoping this recursion gives us that  $v_n^{k+1}$ is described by a degree-$(k-1)$ polynomial, completing the inductive step. 
\end{proof}

\begin{corollary}
    The explicit formula for $b_n^d$ is a degree-$(d-4)$ polynomial, and the explicit formula for $w_n^d$ is a degree-$(d-2)$ polynomial.
\end{corollary}

Theorem \ref{thm:polynomials} and the recurrences for $v_n^d$, $b_n^d$, and $w_n^d$ allow us to compute explicit polynomial formulas for $v_n^d$ using finite differences.

\subsection{Number of vertices of degree $2n-a$ in $\OFG(M_{2,n})$}

We have proven that the number of vertices of a particular degree $d$ is described by a degree-$(d-2)$ polynomial. Interestingly, more polynomial recurrences appear if we express $d$ in terms of $n$ instead of as a fixed natural number. For instance, we can examine the number of vertices in the $n$th generation of $\chi_D$ labeled $2n-a$ for all $a \in \N_0$, where $\N_0=\N \cup \{0\}$. This considers vertices of degree $a$ less than the maximum degree $2n$ of $\OFG(M_{2,n})$. Similar to how we collected sequences to construct Table \ref{tab:degdistribution}, we now collect sequences $(v_n^{2n-a})_n$, this time excluding the edge case where there are exactly $4$ vertices labeled with $2n-a$ (i.e. $2n-a=2$). We then compute the differences between consecutive terms of the sequences to look for polynomial patterns. This gives us the following new table of sequences:

\begin{table}
\centering
\begin{tabular}{c|c|c|c|c}
\textbf{$a$} & \textbf{$(v_n^{2n-a})$} & \textbf{$\Delta(v_n^{2n-a})$} & \textbf{$\Delta^2(v_n^{2n-a})$} & \textbf{$\Delta^3(v_n^{2n-a})$} \\ \hline
\textbf{$0$} & $(2,2,2,2,2,2,2,2)$ & $(0,0,0,0,0,0,0)$ &  &  \\ \hline
\textbf{$1$} & $(0,0,0,0,0,0,0,0)$ & $(0,0,0,0,0,0,0)$ & & \\ \hline
\textbf{$2$} & $(8,12,16,20,24,28,32)$ & $(4,4,4,4,4,4)$ & $(0,0,0,0,0,0)$ &  \\ \hline
\textbf{$3$} & $(4,8,12,16,20,24,28)$ & $(4,4,4,4,4,4)$ & $(0,0,0,0,0)$ &  \\ \hline
\textbf{$4$} & $(20,44,76,116,164,220)$ & (24,32,40,48,56) & $(8,8,8,8,8)$ & $(0,0,0,0)$ \\ \hline
\textbf{$5$} & $(8,36,80,140,216,308)$ & $(28,44,60,76,92)$ & $(16,16,16,16)$ & $(0,0,0)$ \\
\end{tabular}

\caption{Sequences of $(v_n^{2n-a})_n$ for $0 \le a \le 5$ and $2 \le n \le 9$. The code computed all sequences up to $a=16$, but only confirmed polynomial patterns up to $a=6$. The data suggests that for $a \in \N$, $(v_n^{2n-a})$ is described by a degree-$\left\lfloor \frac{a}{2} \right\rfloor$ polynomial.}
\label{tab:reversepolynomialdata}
\end{table}

The data in Table \ref{tab:reversepolynomialdata} suggests the following Theorem, which can be proven using Proposition~\ref{prop:recurrences} and induction.

\begin{theorem} \label{thm:reversepolynomials}
The explicit formula for $v_n^{2n-a}$ where $n \ge \left\lceil \frac{a}{2} \right\rceil+1$ is a polynomial of $n$ with degree $\left\lfloor \frac{a}{2} \right\rfloor$ for all $a \in \N_0$.
\end{theorem}

The techniques we used throughout this section to describe the degree sequences of $\OFG(M_{2,n})$ are clearly specific to the $2\times n$ case. Expanding them to beyond $m = 2$ rows of the Miura-ori would require different methods, and it is not clear if the closed formulas for the number of vertices of a fixed degree $d$ will continue to be polynomial.


\section{Diameter of \texorpdfstring{$\OFG(M_{2,n})$}{OFG(M2n)}}\label{subsec:diameter}
In Lemma \ref{thm:oppositeexist} we showed that there exist exactly four vertices of degree 2 in the origami flip graph of $M_{2,n}$. We now build on this result to explicitly find the diameter of $\OFG(M_{2,n})$. Define two MV assignments with the same flippable faces to be \textit{opposite} if they have opposite MV parities for all their edges. We will find the length of a shortest path between two opposite degree-2 vertices and then show that this is, in fact, the diameter of $\OFG(M_{2,n})$. 

\begin{figure}
\centering
\resizebox{0.7\textwidth}{!}{%
\begin{circuitikz}
\tikzstyle{every node}=[font=\Huge]';
\draw [line width=1.5pt, short] (3.75,20) -- (2.5,17.5);
\draw [line width=1.5pt, short] (2.5,17.5) -- (3.75,15);
\draw [line width=1.5pt, short] (6.25,20) -- (5,17.5);
\draw [line width=1.5pt, dashed] (5,17.5) -- (6.25,15);
\draw [line width=1.5pt, short] (3.75,20) -- (11.25,20);
\draw [line width=1.5pt, short] (11.25,20) -- (10,17.5);
\draw [line width=1.5pt, short] (10,17.5) -- (11.25,15);
\draw [line width=1.5pt, short] (11.25,15) -- (3.75,15);
\draw [line width=1.5pt, short] (8.75,20) -- (7.5,17.5);
\draw [line width=1.5pt, dashed] (7.5,17.5) -- (10,17.5);
\draw [line width=1.5pt, dashed] (2.5,17.5) -- (5,17.5);
\draw [line width=1.5pt, dashed] (7.5,17.5) -- (5,17.5);
\draw [line width=1.5pt, dashed] (7.5,17.5) -- (8.75,15);

\node [font=\Huge, color={rgb,255:red,255; green,0; blue,0}] at (9.5,18.75) {flip};

\draw [line width=1.5pt]  [-{Stealth[length=4mm, width=4mm]}] (10.5,17.5) -- (12,17.5);

\draw [line width=1.5pt, short] (13.75,20) -- (12.5,17.5);
\draw [line width=1.5pt, short] (12.5,17.5) -- (13.75,15);
\draw [line width=1.5pt, short] (16.25,20) -- (15,17.5);
\draw [line width=1.5pt, dashed] (15,17.5) -- (16.25,15);
\draw [line width=1.5pt, short] (13.75,20) -- (21.25,20);
\draw [line width=1.5pt, short] (21.25,20) -- (20,17.5);
\draw [line width=1.5pt, short] (20,17.5) -- (21.25,15);
\draw [line width=1.5pt, short] (21.25,15) -- (13.75,15);
\draw [line width=1.5pt, dashed] (18.75,20) -- (17.5,17.5);
\draw [line width=1.5pt, short] (17.5,17.5) -- (20,17.5);
\draw [line width=1.5pt, dashed] (12.5,17.5) -- (15,17.5);
\draw [line width=1.5pt, dashed] (17.5,17.5) -- (15,17.5);
\draw [line width=1.5pt, dashed] (17.5,17.5) -- (18.75,15);

\node [font=\Huge, color={rgb,255:red,255; green,0; blue,0}] at (14.5,16.25) {flip};

\draw [line width=1.5pt] [-{Stealth[length=4mm, width=4mm]}] (20.5,17.5) -- (22,17.5);

\draw [line width=1.5pt, short] (23.75,20) -- (22.5,17.5);
\draw [line width=1.5pt, short] (22.5,17.5) -- (23.75,15);
\draw [line width=1.5pt, short] (26.25,20) -- (25,17.5);
\draw [line width=1.5pt, short] (25,17.5) -- (26.25,15);
\draw [line width=1.5pt, short] (23.75,20) -- (31.25,20);
\draw [line width=1.5pt, short] (31.25,20) -- (30,17.5);
\draw [line width=1.5pt, short] (30,17.5) -- (31.25,15);
\draw [line width=1.5pt, short] (31.25,15) -- (23.75,15);
\draw [line width=1.5pt, dashed] (28.75,20) -- (27.5,17.5);
\draw [line width=1.5pt, short] (27.5,17.5) -- (30,17.5);
\draw [line width=1.5pt, short] (22.5,17.5) -- (25,17.5);
\draw [line width=1.5pt, dashed] (27.5,17.5) -- (25,17.5);
\draw [line width=1.5pt, dashed] (27.5,17.5) -- (28.75,15);

\node [font=\Huge, color={rgb,255:red,255; green,0; blue,0}] at (27,16.25) {flip};

\draw [line width=1.5pt] (30.5,14.25)  arc
    [  start angle=-90,
        end angle=90,
        radius=1.5cm    ] ;
\draw [line width=1.5pt] (30.5,14.25) -- (2.25,14.25);
\draw [line width=1.5pt] [-{Stealth[length=4mm, width=4mm]}]  (2.25,14.25) arc (90:270:1.5) ;

\draw [line width=1.5pt, short] (3.75,13.75) -- (2.5,11.25);
\draw [line width=1.5pt, short] (2.5,11.25) -- (3.75,8.75);
\draw [line width=1.5pt, short] (6.25,13.75) -- (5,11.25);
\draw [line width=1.5pt, dashed] (5,11.25) -- (6.25,8.75);
\draw [line width=1.5pt, short] (3.75,13.75) -- (11.25,13.75);
\draw [line width=1.5pt, short] (11.25,13.75) -- (10,11.25);
\draw [line width=1.5pt, short] (10,11.25) -- (11.25,8.75);
\draw [line width=1.5pt, short] (11.25,8.75) -- (3.75,8.75);
\draw [line width=1.5pt, dashed] (8.75,13.75) -- (7.5,11.25);
\draw [line width=1.5pt, short] (7.5,11.25) -- (10,11.25);
\draw [line width=1.5pt, short] (2.5,11.25) -- (5,11.25);
\draw [line width=1.5pt, short] (7.5,11.25) -- (5,11.25);
\draw [line width=1.5pt, short] (7.5,11.25) -- (8.75,8.75);

\node [font=\Huge, color={rgb,255:red,255; green,0; blue,0}] at (4.5,12.5) {flip};

\draw [line width=1.5pt]  [-{Stealth[length=4mm, width=4mm]}] (10.5,11.25) -- (12,11.25);

\draw [line width=1.5pt, short] (13.75,13.75) -- (12.5,11.25);
\draw [line width=1.5pt, short] (12.5,11.25) -- (13.75,8.75);
\draw [line width=1.5pt, dashed] (16.25,13.75) -- (15,11.25);
\draw [line width=1.5pt, dashed] (15,11.25) -- (16.25,8.75);
\draw [line width=1.5pt, short] (13.75,13.75) -- (21.25,13.75);
\draw [line width=1.5pt, short] (21.25,13.75) -- (20,11.25);
\draw [line width=1.5pt, short] (20,11.25) -- (21.25,8.75);
\draw [line width=1.5pt, short] (21.25,8.75) -- (13.75,8.75);
\draw [line width=1.5pt, dashed] (18.75,13.75) -- (17.5,11.25);
\draw [line width=1.5pt, short] (17.5,11.25) -- (20,11.25);
\draw [line width=1.5pt, dashed] (12.5,11.25) -- (15,11.25);
\draw [line width=1.5pt, short] (17.5,11.25) -- (15,11.25);
\draw [line width=1.5pt, short] (17.5,11.25) -- (18.75,8.75);

\node [font=\Huge, color={rgb,255:red,255; green,0; blue,0}] at (14.5,10) {flip};

\draw [line width=1.5pt] [-{Stealth[length=4mm, width=4mm]}] (20.5,11.25) -- (22,11.25);

\draw [line width=1.5pt, short] (23.75,13.75) -- (22.5,11.25);
\draw [line width=1.5pt, short] (22.5,11.25) -- (23.75,8.75);
\draw [line width=1.5pt, dashed] (26.25,13.75) -- (25,11.25);
\draw [line width=1.5pt, short] (25,11.25) -- (26.25,8.75);
\draw [line width=1.5pt, short] (23.75,13.75) -- (31.25,13.75);
\draw [line width=1.5pt, short] (31.25,13.75) -- (30,11.25);
\draw [line width=1.5pt, short] (30,11.25) -- (31.25,8.75);
\draw [line width=1.5pt, short] (31.25,8.75) -- (23.75,8.75);
\draw [line width=1.5pt, dashed] (28.75,13.75) -- (27.5,11.25);
\draw [line width=1.5pt, short] (27.5,11.25) -- (30,11.25);
\draw [line width=1.5pt, short] (22.5,11.25) -- (25,11.25);
\draw [line width=1.5pt, short] (27.5,11.25) -- (25,11.25);
\draw [line width=1.5pt, short] (27.5,11.25) -- (28.75,8.75);
\end{circuitikz}
}%
\caption{A shortest sequence of face flips on $M_{2,3}$ between the MV assignments of two degree-2 vertices in $\OFG(M_{2,3})$, showing they have distance 5 in this graph.}\label{example2x3}
\end{figure}
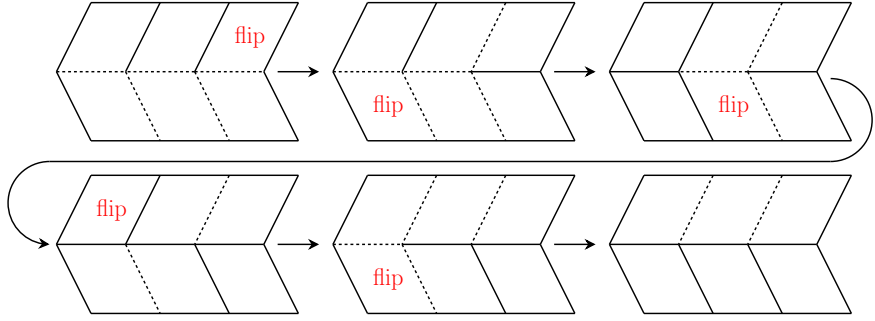

This task is not as easy as one might think. For example, flipping from one degree-2 vertex in $\OFG(M_{2,3})$ to its opposite cannot be done simply by flipping three of the six faces of $M_{2,3}$ to ensure that each crease changes its MV parity, as this would cause a violation of  the Big Little Big Theorem along the way. Figure~\ref{example2x3} shows one way to do it in five flips, which turns out to be the fewest possible. Since these sequences of flips must be done carefully, we turn to results on 3-colorings and use the coloring-Miura bijection  from Section~\ref{sec:bg}.

 In particular, we use techniques and results from \cite{Cereceda1,3colorReconfig} on the shortest paths between 3-colorings in the \textit{3-coloring reconfiguration graph} $R_3(G)$, which is the graph whose vertices are the proper 3-colorings of $G$, and two vertices are adjacent if and only if the two colorings differ at only one vertex of $G$. Define a \textit{$s_0\to s_\ell$ recoloring} as a sequence $R=s_0,...,s_\ell$ of proper colorings of $G$, where $s_{i-1}$ and $s_i$, $1\leq i\leq \ell$, differ in color at exactly one vertex; this corresponds to a path of length $\ell$ in $R_3(G)$. Since changing a vertex color on $M^*_{m,n}$ is equivalent to flipping a face in $M_{m,n}$, a recoloring on $M^*_{m,n}$ is equivalent to a sequence of face flips on $M_{m,n}$. Thus, if we want to find the shortest path between two MV assignments in $\OFG(M_{2,n})$, we look at the 3-colorings in each of their corresponding equivalence class. We find the shortest recoloring between any pair of 3-colorings, one from each equivalence class, and take a pair with the minimum recoloring length, as it corresponds to the shortest distance between the two MV assignments. 
First, we define a few concepts regarding 3-colorings and 3-coloring reconfiguration graphs. 

Suppose $\gamma$ is a proper 3-coloring of a graph $G$ and $u_i,u_j$ are two adjacent vertices in $G$. The weight $w(\gamma,\overrightarrow{u_iu_j})$ of a directed edge from $u_i$ to $u_j$ under $\gamma$ is a value in $\{1,-1\}$ such that $w(\gamma,\overrightarrow{u_iu_j})\equiv \gamma(u_j)-\gamma(u_i) \mod 3$. Let $P_{u,v}$ denote a directed path between two vertices $u,v$ in a 3-colored graph $G$. The \textit{weight} of a path $w(\gamma, P_{u,v})$, is the sum of the weights of the edges in $P_{u,v}$. 

Let $G$ be a graph and $v$ a vertex in $G$. A \textit{simple cycle} is a path of nonzero length from $v$ to $v$ with no repeated edges or vertices beside $v$. 

\begin{lemma}[\cite{cycleLemma}]\label{lem:cycleweight}
    Let $G$ be a graph and let $\gamma,\beta$ be two 3-colorings of $G$. If $\gamma,\beta$ belong to the same component of $R_3(G)$, then, for every simple cycle 
    $K$  of $G$, $w(\gamma,K)=w(\beta,K)$.
\end{lemma}

\begin{lemma}[\cite{tomTessellations}]\label{lem:MiuraOFGconn}
The graph $\OFG(M_{m,n})$, and therefore also $R_3(M_{m,n}^*)$, are connected graphs.
\end{lemma} 

Using these, we can prove the following:

\begin{lemma}\label{lem:weightindependence}
    The weight function on the graph $M_{m,n}^*$ equals zero for all simple cycles and is path independent.

\end{lemma}

\begin{proof}
    Let $K$ be a simple cycle in $M_{m,n}^*$.  By Lemmas \ref{lem:cycleweight} and \ref{lem:MiuraOFGconn}, if we can show that $w(\gamma,K)=0$ for some specific 3-coloring $\gamma$ then the weight of $K$ will equal zero for all 3-colorings. Let $\gamma$ be a 2-coloring of $M_{m,n}^*$ with colors $0$ and $1$ in $\mathbb{Z}_3$, which we know is possible since $M_{m,n}^*$ is bipartite. Then since $K$ is an even cycle, $w(\gamma,K)$ will be the sum of an equal number of $1$s and $-1$s, resulting in zero.

    To prove path independence, we let $P_1$ and $P_2$ be two different directed paths in $M_{m,n}^*$ from vertices $v_{a,b}$ and $v_{c,d}$. We may then use a standard argument of letting $P_2^{-1}$ be the directed path $P_2$ with the orientations reversed, and then note that $P_1\cup P_2^{-1}$ is a union of simple cycles, whose weights must all be zero under 3-colorings,  and duplicate edges in opposite directions.  Thus for any 3-coloring $\gamma$ we have $w(\gamma,P_1)+w(\gamma,P_2^{-1}) = w(\gamma, P_1\cup P_2^{-1}) = 0$, which implies $w(\gamma,P_1) = w(\gamma, P_2)$.
\end{proof}

Lemma~\ref{lem:weightindependence} makes the following  well-defined: Let $\gamma,\beta$ be two 3-colorings of $M_{m,n}^*$. The \textit{relative height} $h_{\gamma,u}(\beta,v)$ of $\beta,v$ with respect to $\gamma,u$ is
$$h_{\gamma,u}(\beta,v)=w(\beta, P_{u,v})-w(\gamma,P_{u,v}),$$
where $P_{u,v}$ is a directed path in $M_{m,n}^*$ from $u$ to $v$. Intuitively, $h_{\gamma,u}(\beta,v)$ is the ``color" difference of a vertex $v$ between $\beta$ and $\gamma$ relative to a vertex $u$. \\ 

We define the \textit{absolute height} of vertex $u$, $H_u^R(s_i)$ with respect to a recoloring $R=s_0, \ldots, s_\ell$ as follows: 
\begin{align*}
    H_u^R(s_0)&=0\\
    H_u^R(s_i)&=\begin{cases}
    H_u^R(s_{i-1}), & \mbox{if }s_i(u)=s_{i-1}(u);\\
    H_u^R(s_{i-1})+2, & \mbox{if }s_i(u)\equiv s_{i-1}(u)-1\mod 3\\
    H_u^R(s_{i-1})-2, & \mbox{if }s_i(u)\equiv s_{i-1}(u)+1\mod 3.
\end{cases}
\end{align*}
Intuitively, if we let $\gamma, s_1,...,s_{\ell-1}, \beta$ be our recoloring from $\gamma \to \beta$, $H_u^R(\beta)$ is the specific color permutation of $\beta$ relative to $\gamma$ at $u$. While $h_{\gamma,u}(\beta,v)$ describes the color relationship between vertices $v$ and $u$ in $\beta$, $H_u^R(\beta)$ describes how the color of $u$ changes. Thus, the $h_{\gamma,v_{1,1}}(\beta,v)$ values describe the equivalence class of 3-colorings $\beta$ occupies, while the $H_{v_{1,1}}^R(\beta)$ value describes if the recoloring $R$ changes the equivalence class we are in.

We have the following lemma from \cite{3colorReconfig}. 

\begin{lemma}[\cite{3colorReconfig}]\label{lem:3colshortestpath}
    Let $G=(V,E)$ be a connected graph with given vertex $u\in V$ and let $\gamma,\beta$ be two 3 colorings of $G$ Let $R=\gamma,s_1,...,s_{\ell-1},\beta$ be a $\gamma\to\beta$ recoloring of length $\ell$. Then 
    $$\ell \geq \frac{1}{2}\sum_{v\in V} |H_u^R(\beta)+h_{\gamma,u}(\beta,v)|$$
    In the case where $\gamma$ and $\beta$ are in the same component of $R_3(G)$, then there exists a path that achieves this lower bound, which equals the distance between $\gamma$ and $\beta$ in $R_3(G)$.
\end{lemma}

We apply this lemma to bound the distance (shortest past length) between any two MV assignments in the general origami flip graph $\OFG(M_{m,n})\cong R_3(M_{m,n}^*)$, and then use this to find the exact distance for $M_{2,n}$. First, we calculate the distance between opposite degree-2 vertices in $\OFG(M_{m,n})$.

\begin{theorem}\label{thm:distbound}
    The length $L$ of the shortest path between two opposite degree-2 vertices in the origami flip graph $\OFG(M_{m,n})$, where $m\leq n$, is
    $$ \left\{
    \begin{array}{cl}
    \frac{mn^2}{4}+\frac{m^3}{12}-\frac{m}{12} & \mbox{if }n-m\mbox{ is odd, and}\\
    \frac{m n^2}{4}+\frac{m^3}{12}-\frac{m}{3} &  \mbox{otherwise.}
    \end{array}
    \right.$$ 
    If $m>n$ then the above formulas hold with $m$ and $n$ reversed.
\end{theorem}
\begin{proof}
We show this distance formula holds for the case where the two flippable faces are $\alpha_{1,n}$ and $\alpha_{m,1}$. The other degree-2 case is a mirror-image of this, with flippable faces $\alpha_{1,1}$ and $\alpha_{m,n}$, so it has the same shortest path length by symmetry. 

By Lemma \ref{lem:3colshortestpath}, we know $\ell \geq \frac{1}{2}\sum_{v\in V(M^*_{m,n})} |H_{v_{1,1}}^R(\beta)+h_{\gamma,v_{1,1}}(\beta,v)|$, so if we can express our problem in terms of three colorings, we’ll have a method for calculating the shortest distance.

Let $\mu$ denote the MV assignment of $M_{m,n}$ where creases $\{\alpha_{i,j},\alpha_{i,j+1}\}$ for $1\leq i\leq m-1$ odd and any $1\leq j\leq n-1$ are mountain and all other creases are valley.
Let $\nu$ be the opposite MV assignment of $\mu$. Recall from Theorem~\ref{thm:4verticesdeg2} and Figure~\ref{fig:deg2vertices} that these are the two MV assignments where $\alpha_{1,n}$ and $\alpha_{m,1}$ are the only flippable faces. Both $\mu$ and $\nu$ have three corresponding 3-colorings of the $m\times n$ grid graph, depending on which color we decide to fix to vertex $v_{1,1}$. Let $\gamma_1,\gamma_2,\gamma_3$ denote the three 3-colorings that correspond to $\mu$ and let $\beta_1,\beta_2,\beta_3$ denote the three 3-colorings that correspond to $\nu$. We want to find the shortest recoloring from any $\gamma_i$ to any $\beta_j$. This is the same as finding the shortest recoloring length from $\gamma=\gamma_1$ to any $\beta_j$ due to color permutations. 

The value of $h_{\gamma,v_{1,1}}(\beta_j,v)$ is the same across the different $\beta_j$’s, since the value of edge weights is not dependent on color permutations. Thus, the $\beta_j$ that ends up being the shortest distance from $\gamma$ is determined by the value of $H_{v_{1,1}}^R(\beta_j)$, which accounts for what color the first vertex is fixed as. To find the shortest recoloring from $\gamma$ to $\beta_j$, we calculate $h_{\gamma,v_{1,1}}(\beta,v)$ from $\gamma$ to an arbitrary $\beta=\beta_j$ and choose the $H=H_{v_{1,1}}^R(\beta)$ value that results in the shortest distance. We can choose this $H$ value since the bijection between the Miura-ori and 3-colorings is one between Miura-ori MV assignments and 3-coloring equivalence classes rather than specific 3-colorings. 

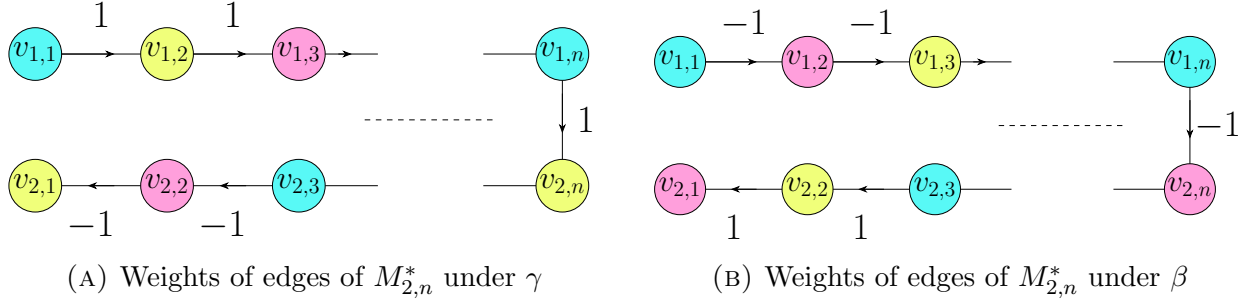
\begin{figure}
    \centering
    \begin{subfigure}[b]{0.48\textwidth}
        \centering
\resizebox{1\textwidth}{!}{%
\begin{circuitikz}
\tikzstyle{every node}=[font=\Large]
\draw [ fill={rgb,255:red,88; green,249; blue,246} ] (6.25,12.25) circle (0.5cm) node {\Large $v_{1,1}$} ;
\draw [ fill={rgb,255:red,240; green,255; blue,122} ] (8.75,12.25) circle (0.5cm) node {\Large $v_{1,2}$} ;
\draw [ fill={rgb,255:red,255; green,159; blue,220} ] (11.25,12.25) circle (0.5cm) node {\Large $v_{1,3}$} ;
\draw [ fill={rgb,255:red,240; green,255; blue,122} ] (6.25,9.75) circle (0.5cm) node {\Large $v_{2,1}$} ;
\draw [ fill={rgb,255:red,255; green,159; blue,220} ] (8.75,9.75) circle (0.5cm) node {\Large $v_{2,2}$} ;
\draw [ fill={rgb,255:red,88; green,249; blue,246} ] (11.25,9.75) circle (0.5cm) node {\Large $v_{2,3}$} ;
\draw [ fill={rgb,255:red,88; green,249; blue,246} ] (6.25,7.25) circle (0.5cm) node {\Large $v_{3,1}$} ;
\draw [ fill={rgb,255:red,240; green,255; blue,122} ] (8.75,7.25) circle (0.5cm) node {\Large $v_{3,2}$} ;
\draw [ fill={rgb,255:red,255; green,159; blue,220} ] (11.25,7.25) circle (0.5cm) node {\Large $v_{3,3}$} ;
\draw [ fill={rgb,255:red,88; green,249; blue,246} ] (16.25,12.25) circle (0.5cm) node {\Large $v_{1,n}$} ;
\draw [ fill={rgb,255:red,255; green,159; blue,220} ] (16.25,9.75) circle (0.5cm) node {\Large $v_{2,n}$} ;
\draw [ fill={rgb,255:red,88; green,249; blue,246} ] (16.25,7.25) circle (0.5cm) node {\Large $v_{3,n}$} ;
\draw [dashed] (12.5,11) -- (15,11);
\draw [dashed] (12.5,8.5) -- (15,8.5);
\draw [short] (6.75,12.25) -- (8.25,12.25);
\draw [short] (9.25,12.25) -- (10.75,12.25);
\draw [short] (16.25,11.75) -- (16.25,10.25);
\draw [short] (15.75,12.25) -- (14.75,12.25);
\draw [short] (15.75,9.75) -- (14.75,9.75);
\draw [short] (15.75,7.25) -- (14.75,7.25);
\draw [short] (11.75,12.25) -- (12.75,12.25);
\draw [short] (11.75,9.75) -- (12.75,9.75);
\draw [short] (10.75,9.75) -- (9.25,9.75);
\draw [short] (8.25,9.75) -- (6.75,9.75);
\draw [->, >=Stealth] (6.75,12.25) -- (7.75,12.25);
\draw [->, >=Stealth] (9.25,12.25) -- (10.25,12.25);
\draw [->, >=Stealth] (12,12.25) -- (12.25,12.25);
\draw [short] (6.75,7.25) -- (8.25,7.25);
\draw [short] (9.25,7.25) -- (10.75,7.25);
\draw [short] (11.75,7.25) -- (12.75,7.25);
\draw [->, >=Stealth] (6.75,7.25) -- (7.75,7.25);
\draw [->, >=Stealth] (9.25,7.25) -- (10.25,7.25);
\draw [->, >=Stealth] (12,7.25) -- (12.25,7.25);
\draw [->, >=Stealth] (16.25,11.5) -- (16.25,10.75);
\draw [->, >=Stealth] (10.25,9.75) -- (9.75,9.75);
\draw [->, >=Stealth] (7.75,9.75) -- (7.25,9.75);
\draw [short] (6.25,9.25) -- (6.25,7.75);
\draw [->, >=Stealth] (6.25,9.25) -- (6.25,8.25);
\draw [->, >=Stealth] (12.25,9.75) -- (12,9.75) ;
\node [font=\LARGE] at (7.5,13) {$1$};
\node [font=\LARGE] at (10,13) {$1$};
\node [font=\LARGE] at (16.7,11) {$-1$};
\node [font=\LARGE] at (9.8,10.5) {$-1$};
\node [font=\LARGE] at (7.3,10.5) {$-1$};
\node [font=\LARGE] at (7.5,8) {$1$};
\node [font=\LARGE] at (10,8) {$1$};
\node [font=\LARGE] at (6.75,8.5) {$-1$};
\end{circuitikz}}%
\caption{Weights of edges of $M_{3,n}^*$ under $\gamma$}
    \end{subfigure}
    \hfill
    \begin{subfigure}[b]{0.48\textwidth}
    \centering
    \resizebox{\textwidth}{!}{%
    \begin{circuitikz}
\draw [ fill={rgb,255:red,88; green,249; blue,246} ] (6.25,5.25) circle (0.5cm) node {\Large $v_{1,1}$} ;
\draw [ fill={rgb,255:red,255; green,159; blue,220} ] (8.75,5.25) circle (0.5cm) node {\Large $v_{1,2}$} ;
\draw [ fill={rgb,255:red,240; green,255; blue,122} ] (11.25,5.25) circle (0.5cm) node {\Large $v_{1,3}$} ;
\draw [ fill={rgb,255:red,255; green,159; blue,220} ] (6.25,2.75) circle (0.5cm) node {\Large $v_{2,1}$} ;
\draw [ fill={rgb,255:red,240; green,255; blue,122} ] (8.75,2.75) circle (0.5cm) node {\Large $v_{2,2}$} ;
\draw [ fill={rgb,255:red,88; green,249; blue,246} ] (11.25,2.75) circle (0.5cm) node {\Large $v_{2,3}$} ;
\draw [ fill={rgb,255:red,88; green,249; blue,246} ] (6.25,0.25) circle (0.5cm) node {\Large $v_{3,1}$} ;
\draw [ fill= {rgb,255:red,255; green,159; blue,220}  ] (8.75,0.25) circle (0.5cm) node {\Large $v_{3,2}$} ;
\draw [ fill={rgb,255:red,240; green,255; blue,122} ] (11.25,0.25) circle (0.5cm) node {\Large $v_{3,3}$} ;
\draw [ fill={rgb,255:red,88; green,249; blue,246} ] (16.25,5.25) circle (0.5cm) node {\Large $v_{1,n}$} ;
\draw [ fill={rgb,255:red,240; green,255; blue,122} ] (16.25,2.75) circle (0.5cm) node {\Large $v_{2,n}$} ;
\draw [ fill={rgb,255:red,88; green,249; blue,246} ] (16.25,0.25) circle (0.5cm) node {\Large $v_{3,n}$} ;
\draw [dashed] (12.5,4) -- (15,4);
\draw [dashed] (12.5,1.5) -- (15,1.5);
\draw [short] (6.75,5.25) -- (8.25,5.25);
\draw [short] (9.25,5.25) -- (10.75,5.25);
\draw [short] (16.25,4.75) -- (16.25,3.25);
\draw [short] (15.75,5.25) -- (14.75,5.25);
\draw [short] (15.75,2.75) -- (14.75,2.75);
\draw [short] (15.75,0.25) -- (14.75,0.25);
\draw [short] (11.75,5.25) -- (12.75,5.25);
\draw [short] (11.75,2.75) -- (12.75,2.75);
\draw [->, >=Stealth] (12.25,2.75) -- (12,2.75) ;
\draw [short] (10.75,2.75) -- (9.25,2.75);
\draw [short] (8.25,2.75) -- (6.75,2.75);
\draw [->, >=Stealth] (6.75,5.25) -- (7.75,5.25);
\draw [->, >=Stealth] (9.25,5.25) -- (10.25,5.25);
\draw [->, >=Stealth] (12,5.25) -- (12.25,5.25);
\draw [->, >=Stealth] (16.25,4.5) -- (16.25,3.75);
\draw [->, >=Stealth] (10.25,2.75) -- (9.75,2.75);
\draw [->, >=Stealth] (7.75,2.75) -- (7.25,2.75);
\draw [short] (6.25,2.25) -- (6.25,0.75);
\draw [->, >=Stealth] (6.25,2.25) -- (6.25,1.25);
\draw [short] (6.75,0.25) -- (8.25,0.25);
\draw [short] (9.25,0.25) -- (10.75,0.25);
\draw [->, >=Stealth] (6.75,0.25) -- (7.75,0.25);
\draw [->, >=Stealth] (9.25,0.25) -- (10.25,0.25);
\draw [short] (11.75,0.25) -- (12.75,0.25);
\draw [->, >=Stealth] (12,0.25) -- (12.25,0.25);
\node [font=\LARGE] at (7.5,6) {$-1$};
\node [font=\LARGE] at (10,6) {$-1$};
\node [font=\LARGE] at (16.8,4) {$1$};
\node [font=\LARGE] at (9.8,3.5) {$1$};
\node [font=\LARGE] at (7.3,3.5) {$1$};
\node [font=\LARGE] at (6.75,1.5) {$1$};
\node [font=\LARGE] at (7.5,1) {$-1$};
\node [font=\LARGE] at (10,1) {$-1$};
   \end{circuitikz}}%
    \caption{Weights of edges of $M_{3,n}^*$ under $\beta$}
    \end{subfigure}
    \caption{3-colorings $\gamma$, $\beta$, corresponding to $\mu$ and $\nu$}
    \label{fig:edgeweights}
\end{figure}

We first calculate $h_{\gamma,v_{1,1}}(\beta,v)$. We construct our path in a boustrophedon pattern in the same way as the construction of the bijection between Miura-ori MV assignments and 3-colored grid graphs (see Figure \ref{fig:edgeweights}).  
By construction and our choice of $\gamma$, 
$$\begin{matrix}\begin{cases}
w(\gamma,\overrightarrow{v_{i,j}v_{i,j+1}})= 1, & 1\leq j<n\\
w(\gamma,\overrightarrow{v_{i,n}v_{i+1,n}})=-1\\
w(\gamma,\overrightarrow{v_{i+1,j}v_{i+1,j-1}})=-1, & 1<j\leq n\\
w(\gamma,\overrightarrow{v_{i+1,1}v_{i+2,1}})=-1
\end{cases} &\text{    } &\begin{cases}
w(\beta,\overrightarrow{v_{i,j}v_{i,j+1}})= -1, & 1\leq j<n\\
w(\beta,\overrightarrow{v_{i,n}v_{i+1,n}})=1\\
w(\beta,\overrightarrow{v_{i+1,j}v_{i+1,j-1}})=1, & 1<j\leq n\\
w(\gamma,\overrightarrow{v_{i+1,1}v_{i+2,1}})=1
\end{cases}\end{matrix}$$
where $i$ is odd ($1\leq i< m-1$). Since $h_{\gamma,v_{1,1}}(\beta,v_{i,j})=w(\beta, P_{v_{1,1},v_{i,j}})-w(\gamma,P_{v_{1,1},v_{i,j}})$ where we take $P_{v_{1,1},v_{i,j}}$ to be our directed boustrophedon path on $M^*_{m,n}$ (as seen in Figure~\ref{fig:3colorbiject}), we have 
\begin{eqnarray}\label{eq:h}
h_{\gamma,v_{1,1}}(\beta,v_{i,j}) & = & 
 \begin{cases}
2(i-1)-2(j-1) = 2i-2j& i\mbox{ odd}, 1\leq j\leq n\\
2i-2(n+1)+2(n-(j-1))=2i-2j & i\mbox{ even}, 1\leq j\leq n.
\end{cases}
\end{eqnarray}
But this just means $h_{\gamma,v_{1,1}}(\beta,v_{i,j})=2i-2j$ in both the $i$ even and $i$ odd cases.
We then calculate the minimum distance $\ell$ using Lemma~\ref{lem:3colshortestpath} and dividing into three cases: $H\leq -2m, -2m<H\leq 2n$, and $H>2n$. In each case, we will need to see when $H+h_{\gamma,v_{1,1}}(\beta,v)$ is positive or negative in order to simplify the absolute value in the summation of Lemma~\ref{lem:3colshortestpath}. 

\textit{Case 1:} $H\leq -2m$.  

Since $1\leq i\leq m$, we have $H-2i-2j\leq H+2m-2j < 0$ in this case. Then
 Lemma~\ref{lem:3colshortestpath} gives us
\begin{align*}
    \ell&\geq \frac{1}{2}\sum_{v\in V(M^*_{m,n})} |H+h_{\gamma,v_{1,1}}(\beta,v)|\\
    & = \frac{1}{2}\sum_{i=1}^m\sum_{j=1}^n -(H+2i-2j)
    =\frac{mn^2}{2}-\frac{m^2 n}{2}-\frac{1}{2}Hmn.
\end{align*}
This reaches its minimum when $H=-2m$, giving us $\ell \geq m n^2/2 + m^2 n/2$.

\vspace{.1in}

\textit{Case 2:} $-2m < H\leq 2n$.

If $1\leq j\leq H/2 + i$ we have
$$H+h_{\gamma,v_{1,1}}(\beta,v_{i,j}) = H+2i-2j \geq H+2i-H-2i=0,$$
whereas if $H/2+i< j\leq n$ we have $H+h_{\gamma,v_{1,1}}(\beta,v_{i,j})$ is negative.  However, if $H<0$ then the $1\leq j\leq H/2 + i$ case might not occur for some values of $i$, so let us assume for now that $H\geq 0$. 
Then Lemma~\ref{lem:3colshortestpath} gives   
\begin{align*}
    \ell&\geq \frac{1}{2}\sum_{v\in V(M^*_{m,n})} |H+h_{\gamma,v_{1,1}}(\beta,v)|\\
    & = \frac{1}{2}\sum_{i=1}^m
    \left(
    \sum_{j=1}^{H/2+i}(H+2i-2j)+
    \sum_{j=H/2+i+1}^n -(H+2i-2j)  
    \right).
\end{align*}
This simplifies to $m/12 (6n^2-6mn-6Hn+4m^2+6Hm+3H^2-4)$, a concave-up quadratic in $H$ which has a minimum at $H=n-m$ and then gives a minimum bound of $\ell\geq m n^2/4+m^3/12-m/3$.  Since we are assuming $H\geq 0$, this formula will hold when $m\leq n$.

However, $H$ is always even (as can be seen in the definition of $H_u^R(s_i)$), so if $n-m$ is odd then $H$ cannot achieve this minimum value.  Instead, the minimum for $H$ will be at $H=m-n+1$ and $m-n-1$, both of which will give the same bound by the symmetry of the quadratic in $H$ and which gives the bound $\ell\geq m n^2/4 + m^3/12 - m/12$.

For the $H<0$ case, we use the fact that the grid graphs $M_{m,n}^*$ and $M_{n,m}^*$ are isomorphic and ``transposing" the graph in this way only reverses the sign of $h$. Therefore the $H<0$ case will give the same bounds but with $m$ and $n$ reversed.

\vspace{.1in}

\textit{Case 3:} $H > 2(n-1)$.

This case turns out to be the same as Case 1. Since $H>2n$ we have $H+2i-2j>0$.  Lemma~\ref{lem:3colshortestpath} then says
\begin{align*}
    \ell&\geq \frac{1}{2}\sum_{v\in V(M^*_{m,n})} |H+h_{\gamma,v_{1,1}}(\beta,v)|\\
    & = \frac{1}{2}\sum_{i=1}^m\sum_{j=1}^n (H+2i-2j) =-\frac{mn^2}{2}+\frac{m^2 n}{2}+\frac{1}{2}Hmn.
\end{align*}
This achieves its minimum when $H=2n$ giving us $\ell \geq mn^2/2 + m^2n/2$, which is the same bound as Case 1.

It is not hard to show that for $m\leq n$,  the Case 2 bounds of $mn^2/4+m^3/12-m/3$ and $mn^2/4+m^3/12-m/12$ are both less than $mn^2/2+m^2n/2$ from Cases 1 and 3. This gives for formula for $L$ as stated in the Theorem.  

Finally, by Lemma~\ref{lem:MiuraOFGconn} we know that the Miura-ori flip graph is connected, so by Lemma \ref{lem:3colshortestpath} there exists a shortest path in $R_3(M_{m,n}^*) \cong \OFG(M_{m,n})$ between opposite degree $2$ vertices that is exactly length $L$.

\end{proof}

Therefore, the shortest distance between the opposite degree-2 vertices in $\OFG(M_{2,n})$ is $\lceil \frac{n^2}{2}\rceil$.  We will now show that the shortest distance between any two vertices in $\OFG(M_{2,n})$ is never longer than this distance (and then discuss why our approach will not work for $M_{m,n}$ with $m>2$).

Referring back to Lemma \ref{lem:3colshortestpath}, we want to know at which $H=H_{v_{1,1}}^R(\beta)$ value $\sum_{v\in V(M^*_{2,n})} |H+h_{\gamma,v_{1,1}}(\beta,v)|$ is minimized, where $\gamma$ and $\beta$ are corresponding 3-colorings of two arbitrary MV assignments on $M_{2,n}$, as that would determine the shortest path length between any two MV assignments. Since $H$ is constant in this summation, we turn our attention to the $h_{\gamma,v_{1,1}}(\beta,v)$ term. 

\begin{lemma}\label{lem:2diff}
For fixed 3-colorings $\gamma$ and $\beta$ of $M^*_{2,n}$, the value of $h_{\gamma,v_{1,1}}(\beta,v)$ changes by 0, $-2$, or 2 as we move from a vertex $v$ in $M^*_{2,n}$ to one of its neighbors.
\end{lemma}

\begin{proof}
    This is a simple consequence of the definition $h_{\gamma,v_{1,1}}(\beta,v)=w(\beta,P_{v_{1,1}, v})-w(\gamma,P_{v_{1,1}, v})$. Since our 3-colorings are proper, as we move from $v$ to one of its neighbors, $w(\beta,P_{v_{1,1}, v})$ and $w(\gamma,P_{v_{1,1}, v})$ will both change by $\pm 1$. If both change by the same amount, $h_{\gamma,v_{1,1}}(\beta,v)$ will not change. If they change by different amounts, $h_{\gamma,v_{1,1}}(\beta,v)$ will change by 2.
\end{proof}

\begin{theorem}\label{thm:shortestpath}
    In $\OFG(M_{2,n})$, the two vertices with a maximal shortest path between them are two opposite vertices of degree $2$. In other words, if we have any two MV assignments $\mu_1,\mu_2$ of $M_{2,n}$, then there exist corresponding 3-colorings $\gamma,\beta$ such that   
$$\frac{1}{2}\sum_{v\in V(M^*_{2,n})}|H+h_{\gamma,v_{1,1}}(\beta,v)|\leq \left\lceil \frac{n^2}{2}\right\rceil,$$
where $H$ is chosen such that $\frac{1}{2}\sum_{v\in V(M^*_{2,n})}|H+h_{\gamma,v_{1,1}}(\beta,v)|$ is minimized. 
\end{theorem}

\begin{figure}
    \centering
    \subcaptionbox{Weights of edges of $M_{2,4}^*$ under $\gamma$}[0.3\textwidth]{
    \resizebox{0.3\textwidth}{!}{%
\begin{circuitikz}
\tikzstyle{every node}=[font=\Large]
\draw [ fill={rgb,255:red,88; green,249; blue,246} ] (6.25,5.25) circle (0.5cm) node {\Large $v_{1,1}$} ;
\draw [ fill={rgb,255:red,240; green,255; blue,122} ] (8.75,5.25) circle (0.5cm) node {\Large $v_{1,2}$} ;
\draw [ fill= {rgb,255:red,88; green,249; blue,246}] (11.25,5.25) circle (0.5cm) node {\Large $v_{1,3}$} ;
\draw [ fill={rgb,255:red,255; green,159; blue,220} ] (6.25,2.75) circle (0.5cm) node {\Large $v_{2,1}$} ;
\draw [ fill={rgb,255:red,88; green,249; blue,246} ] (8.75,2.75) circle (0.5cm) node {\Large $v_{2,2}$} ;
\draw [ fill={rgb,255:red,255; green,159; blue,220} ] (11.25,2.75) circle (0.5cm) node {\Large $v_{2,3}$} ;
\draw [ fill={rgb,255:red,255; green,159; blue,220} ] (13.75,5.25) circle (0.5cm) node {\Large $v_{1,4}$} ;
\draw [ fill={rgb,255:red,240; green,255; blue,122} ] (13.75,2.75) circle (0.5cm) node {\Large $v_{2,4}$} ;
\draw [short] (6.75,5.25) -- (8.25,5.25);
\draw [short] (9.25,5.25) -- (10.75,5.25);
\draw [short] (13.75,4.75) -- (13.75,3.25);
\draw [short] (11.75,5.25) -- (13.25,5.25);
\draw [short] (11.75,2.75) -- (13.25,2.75);
\draw [->, >=Stealth] (12.5,2.75) -- (12.25,2.75) ;
\draw [short] (10.75,2.75) -- (9.25,2.75);
\draw [short] (8.25,2.75) -- (6.75,2.75);
\draw [->, >=Stealth] (6.75,5.25) -- (7.5,5.25);
\draw [->, >=Stealth] (9.25,5.25) -- (10,5.25);
\draw [->, >=Stealth] (12.25,5.25) -- (12.5,5.25);
\draw [->, >=Stealth] (13.75,4.5) -- (13.75,3.75);
\draw [->, >=Stealth] (10.25,2.75) -- (9.75,2.75);
\draw [->, >=Stealth] (7.75,2.75) -- (7.25,2.75);
\node [font=\LARGE] at (7.4,5.75) {$1$};
\node [font=\LARGE] at (9.9,5.75) {$-1$};
\node [font=\LARGE] at (12.4,5.75) {$-1$};
\node [font=\LARGE] at (14.75,4) {$-1$};
\node [font=\LARGE] at (9.6,3.15) {$1$};
\node [font=\LARGE] at (7.1,3.15) {$-1$};
\node [font=\LARGE] at (12.1,3.15) {$1$};

\draw [line width=1pt, short] (7.25,4) -- (8.5,7);
\draw [line width=1pt, dashed] (9.75,4) -- (11,7);
\draw [line width=1pt, dashed] (12.25,4) -- (13.5,7);
\draw [line width=1pt, dashed] (7.25,4) -- (8.5,1);
\draw [line width=1pt, short] (9.75,4) -- (11,1);
\draw [line width=1pt, short] (12.25,4) -- (13.5,1);
\draw [line width=1pt, dashed] (5.5,4) -- (7.25,4);
\draw [line width=1pt, dashed] (7.25,4) -- (9.75,4);
\draw [line width=1pt, dashed] (9.75,4) -- (12.25,4);
\draw [line width=1pt, dashed] (12.25,4) -- (14.25,4);
\end{circuitikz}}}
\subcaptionbox{Weights of edges of $M^*_{2,4}$ under $\beta$ and the $h_{\gamma,v_{1,1}}({\beta},v_{i,j})$ values (in bold).}[0.3\textwidth]{
\resizebox{0.3\textwidth}{!}{%
\begin{circuitikz}
\draw [ fill={rgb,255:red,88; green,249; blue,246} ] (6.25,5.25) circle (0.5cm) node {\Large $v_{1,1}$} ;
\draw [ fill={rgb,255:red,255; green,159; blue,220} ] (8.75,5.25) circle (0.5cm) node {\Large $v_{1,2}$} ;
\draw [ fill={rgb,255:red,240; green,255; blue,122} ] (11.25,5.25) circle (0.5cm) node {\Large $v_{1,3}$} ;
\draw [ fill={rgb,255:red,255; green,159; blue,220} ] (6.25,2.75) circle (0.5cm) node {\Large $v_{2,1}$} ;
\draw [ fill={rgb,255:red,88; green,249; blue,246} ] (8.75,2.75) circle (0.5cm) node {\Large $v_{2,2}$} ;
\draw [ fill={rgb,255:red,255; green,159; blue,220} ] (11.25,2.75) circle (0.5cm) node {\Large $v_{2,3}$} ;
\draw [ fill={rgb,255:red,255; green,159; blue,220} ] (13.75,5.25) circle (0.5cm) node {\Large $v_{1,4}$} ;
\draw [ fill={rgb,255:red,88; green,249; blue,246} ] (13.75,2.75) circle (0.5cm) node {\Large $v_{2,4}$} ;
\draw [short] (6.75,5.25) -- (8.25,5.25);
\draw [short] (9.25,5.25) -- (10.75,5.25);
\draw [short] (13.75,4.75) -- (13.75,3.25);
\draw [short] (11.75,5.25) -- (13.25,5.25);
\draw [short] (11.75,2.75) -- (13.25,2.75);
\draw [->, >=Stealth] (12.5,2.75) -- (12.25,2.75) ;
\draw [short] (10.75,2.75) -- (9.25,2.75);
\draw [short] (8.25,2.75) -- (6.75,2.75);
\draw [->, >=Stealth] (6.75,5.25) -- (7.5,5.25);
\draw [->, >=Stealth] (9.25,5.25) -- (10,5.25);
\draw [->, >=Stealth] (12.25,5.25) -- (12.5,5.25);
\draw [->, >=Stealth] (13.75,4.5) -- (13.75,3.75);
\draw [->, >=Stealth] (10.25,2.75) -- (9.75,2.75);
\draw [->, >=Stealth] (7.75,2.75) -- (7.25,2.75);
\node [font=\LARGE] at (7.4,5.75) {$-1$};
\node [font=\LARGE] at (9.9,5.75) {$-1$};
\node [font=\LARGE] at (12.4,5.75) {$1$};
\node [font=\LARGE] at (14.75,4) {$1$};
\node [font=\LARGE] at (9.75,3.25) {$1$};
\node [font=\LARGE] at (7.2,3.25) {$-1$};
\node [font=\LARGE] at (12.1,3.25) {$-1$};
\draw [line width=1pt, dashed] (7.25,4) -- (8.5,7);
\draw [line width=1pt, dashed] (9.75,4) -- (11,7);
\draw [line width=1pt, short] (12.25,4) -- (13.5,7);
\draw [line width=1pt, dashed] (7.25,4) -- (8.5,1);
\draw [line width=1pt, short] (9.75,4) -- (11,1);
\draw [line width=1pt, dashed] (12.25,4) -- (13.5,1);
\draw [line width=1pt, dashed] (5.5,4) -- (7.25,4);
\draw [line width=1pt, short] (7.25,4) -- (9.75,4);
\draw [line width=1pt, short] (9.75,4) -- (12.25,4);
\draw [line width=1pt, short] (12.25,4) -- (14.25,4);
\node [font=\LARGE] at (6.25,6.25) {\textbf{0}};
\node [font=\LARGE] at (8.75,6.25) {$\mathbf{-2}$};
\node [font=\LARGE] at (11.25,6.25) {$\mathbf{-2}$};
\node [font=\LARGE] at (13.75,6.25) {$\mathbf{0}$};
\node [font=\LARGE] at (13.75,1.85) {$\mathbf{2}$};
\node [font=\LARGE] at (11.25,1.85) {$\mathbf{0}$};
\node [font=\LARGE] at (8.75,1.85) {$\mathbf{0}$};
\node [font=\LARGE] at (6.25,1.85) {$\mathbf{0}$};
   \end{circuitikz}}}
\subcaptionbox{Graph of the  $ |{h}_{{\gamma},v_{1,1}}({\beta},v_{i,j})|$ values as $v_{i,j}$ travels along the directed cycle $K^*$. Removing the "plateaus" gives a
Dyck path. The area under this represents $\sum |{h}_{{\gamma},v_{1,1}}({\beta},v_{i,j})|$.}[0.3\textwidth]{
\resizebox{0.3\textwidth}{!}{%
    \begin{circuitikz}
     \draw [line width=1pt, short] (0,1) -- (9.5,1);
     \draw [line width=1pt, short] (1,0) -- (1,3);
     \draw [line width=1pt, short] (1,1) -- (2,2) -- (3,2) -- (4,1) -- (5,2) -- (6,1) -- (9,1);
     \draw [line width=1pt, short] (2,.7) -- (2,1.3);
     \draw [line width=1pt, short] (3,.7) -- (3,1.3);
     \draw [line width=1pt, short] (4,.7) -- (4,1.3);
     \draw [line width=1pt, short] (5,.7) -- (5,1.3);
     \draw [line width=1pt, short] (6,.7) -- (6,1.3);
     \draw [line width=1pt, short] (7,.7) -- (7,1.3);
     \draw [line width=1pt, short] (8,.7) -- (8,1.3);
     \draw [line width=1pt, short] (9,.7) -- (9,1.3);
     \draw [line width=1pt, short] (.7,2) -- (1.3,2);

     \node [font=\LARGE] at (10,1) {$v_{i,j}$};
     \node [font=\LARGE] at (2,.5) {$v_{1,2}$};
     \node [font=\LARGE] at (3,.5) {$v_{1,3}$};
     \node [font=\LARGE] at (4,.5) {$v_{1,4}$};
     \node [font=\LARGE] at (5,.5) {$v_{2,4}$};
     \node [font=\LARGE] at (6,.5) {$v_{2,3}$};
     \node [font=\LARGE] at (7,.5) {$v_{2,2}$};
     \node [font=\LARGE] at (8,.5) {$v_{2,1}$};
     \node [font=\LARGE] at (9,.5) {$v_{1,1}$};
     \node [font=\LARGE] at (.5,2) {2};

     \node [font=\LARGE] at (1,3.6) {$|h_{\gamma,v_{1,1}}({\beta},v_{i,j})|$};
     \draw [ fill={rgb,255:red,0; green,0; blue,0} ] (1,1) circle (0.25cm) ; 
     \draw [ fill={rgb,255:red,0; green,0; blue,0} ] (2,2) circle (0.25cm) ; 
     \draw [ fill={rgb,255:red,0; green,0; blue,0} ] (3,2) circle (0.25cm) ; 
     \draw [ fill={rgb,255:red,0; green,0; blue,0} ] (4,1) circle (0.25cm) ; 
     \draw [ fill={rgb,255:red,0; green,0; blue,0} ] (5,2) circle (0.25cm) ; 
     \draw [ fill={rgb,255:red,0; green,0; blue,0} ] (6,1) circle (0.25cm) ; 
     \draw [ fill={rgb,255:red,0; green,0; blue,0} ] (7,1) circle (0.25cm) ; 
     \draw [ fill={rgb,255:red,0; green,0; blue,0} ] (8,1) circle (0.25cm) ; 
     \draw [ fill={rgb,255:red,0; green,0; blue,0} ] (9,1) circle (0.25cm) ; 

\draw [->, >=Stealth] (5,-0.5) -- (5,-1.5);


\begin{scope}[yshift=-5cm]
\draw [line width=1pt, short] (0,1) -- (9.5,1);
     \draw [line width=1pt, short] (1,0) -- (1,3);
     \draw [line width=1pt, short] (1,1) -- (2,2) -- (3,1) -- (4,2) -- (5,1);
     \draw [line width=1pt, short] (2,.7) -- (2,1.3);
     \draw [line width=1pt, short] (3,.7) -- (3,1.3);
     \draw [line width=1pt, short] (4,.7) -- (4,1.3);
     \draw [line width=1pt, short] (5,.7) -- (5,1.3);
     \draw [line width=1pt, short] (6,.7) -- (6,1.3);
     \draw [line width=1pt, short] (7,.7) -- (7,1.3);
     \draw [line width=1pt, short] (8,.7) -- (8,1.3);
     \draw [line width=1pt, short] (9,.7) -- (9,1.3);
     \draw [line width=1pt, short] (.7,2) -- (1.3,2);

     \node [font=\LARGE] at (.5,2) {2};

     \node [font=\LARGE] at (1,3.6) {$|h_{\gamma,v_{1,1}}({\beta},v_{i,j})|$};
     \draw [ fill={rgb,255:red,0; green,0; blue,0} ] (1,1) circle (0.25cm) ; 
     \draw [ fill={rgb,255:red,0; green,0; blue,0} ] (2,2) circle (0.25cm) ; 
     \draw [ fill={rgb,255:red,0; green,0; blue,0} ] (3,1) circle (0.25cm) ; 
     \draw [ fill={rgb,255:red,0; green,0; blue,0} ] (4,2) circle (0.25cm) ; 
     \draw [ fill={rgb,255:red,0; green,0; blue,0} ] (5,1) circle (0.25cm) ; 
\end{scope}

     \end{circuitikz}
     }}%
    \caption{3-colorings of two non opposite, non degree 2 MV assignments $\gamma$, ${\beta}$, and the ${h}_{{\gamma},v_{1,1}}({\beta},v_{i,j})$ values in $M^*_{2,4}$. Notice how the area under the path is smaller than for the case of $\tilde{\gamma}$ and $\tilde{\beta}$.}
    \label{fig:nonmaximalpath}
\end{figure}

\begin{proof}
    Let $\tilde{\gamma},\tilde{\beta}$ be the 3-colorings of $M_{2,n}^*$ corresponding to the MV assignments where only $\alpha_{1,n}$ and $\alpha_{2,1}$ are flippable. We know $\frac{1}{2}\sum_{v\in V(M^*_{2,n})}|H+h_{\tilde{\gamma},v_{1,1}}(\tilde{\beta},v)|=\lceil \frac{n^2}{2}\rceil$. To show this is the diameter, we want to show $\frac{1}{2}\sum_{v\in V(M^*_{2,n})}|H+h_{\gamma,v_{1,1}}(\beta,v)|\leq\frac{1}{2}\sum_{v\in V(M^*_{2,n})}|H+h_{\tilde{\gamma},v_{1,1}}(\tilde{\beta},v)|$, for any two 3-colorings $\gamma,\beta$ that correspond to MV assignments $\mu_1,\mu_2$.


Let $\tilde{h}_{i,j}=h_{\tilde{\gamma},v_{1,1}}(\tilde{\beta},v_{i,j})$ and let $h_{i,j}=h_{\gamma,v_{1,1}}(\beta,v_{i,j})$ where $\gamma$ and $\beta$ are any two 3-colorings of $M^*_{2,n}$.
We know from Lemma~\ref{lem:weightindependence} that the values of $\tilde{h}_{i,j}$ and $h_{i,j}$ are independent of the path taken from $v_{1,1}$ to $v_{i,j}$, so we will assume it is the path from $v_{1,1}$ to $v_{i,j}$ along the directed Hamilton cycle $K^*$  on $M^*_{2,n}$ consisting of edges $\overrightarrow{v_{1,1}v_{1,2}}$,$\overrightarrow{v_{1,2}v_{1,3}}$, $\ldots$, $\overrightarrow{v_{1,n}v_{2,n}}$, $\overrightarrow{v_{2,n}v_{2,n-1}}, \ldots,\allowbreak\overrightarrow{v_{2,2}v_{2,1}},\overrightarrow{v_{2,1}v_{1,1}}$ (which is just our boustropheron path on $M^*_{2,n}$, as in Figure~\ref{fig:3colorbiject}). We know by Lemma~\ref{lem:2diff} that from one vertex to the next along $K^*$ the values of $\tilde{h}_{i,j}$ and $h_{i,j}$ change by 2, 0, or $-2$. As such, the graph of $|h_{i,j}|$ as we travel along $K^*$ starting at $v_{1,1}$ will be like a Dyck path, aside from having ``plateaus" whenever $h_{i,j}$ doesn't change; see Figure~\ref{fig:nonmaximalpath} for an example, as well as how these ``plateaus" could be altered to make a traditional Dyck path.  The sum $\sum |h_{i,j}|$ is an approximation of the area under the corresponding Dyck path, and it is known that a Dyck path's area is maximized when it monotonically increases to a maximum point in the middle of its domain and then monotonically decreases back to zero.\footnote{To see this, consider any other Dyck path, which must have at least one ``valley" point shaped like $\searrow\nearrow$, which could be flipped up to $\nearrow\searrow$ and thereby increase the area under the path.} This is exactly what $|\tilde{h}_{i,j}|$ does along the path $K^*$ (the $M^*_{2,4}$ case is illustrated in Figure~\ref{fig:M24}), proving that 
$\frac{1}{2}\sum_{v\in V(M^*_{2,n})}|H+h_{i,j}|\leq\frac{1}{2}\sum_{v\in V(M^*_{2,n})}|H+\tilde{h}_{i,j}|$ since $H$ is constant. Thus, the maximum distance between two vertices in $\OFG(M_{2,n})$ can be found between the two vertices that represent the MV assignment whose only flippable faces are $\alpha_{1,n}$ and $\alpha_{2,1}$ and its opposite MV assignment.


\begin{figure}
    \centering
    \subcaptionbox{Weights of edges of $M_{2,4}^*$ under $\tilde{\gamma}$}[0.3\textwidth]{
    \resizebox{0.3\textwidth}{!}{%
\begin{circuitikz}
\tikzstyle{every node}=[font=\Large]
\draw [ fill={rgb,255:red,88; green,249; blue,246} ] (6.25,5.25) circle (0.5cm) node {\Large $v_{1,1}$} ;
\draw [ fill={rgb,255:red,240; green,255; blue,122} ] (8.75,5.25) circle (0.5cm) node {\Large $v_{1,2}$} ;
\draw [ fill= {rgb,255:red,255; green,159; blue,220}] (11.25,5.25) circle (0.5cm) node {\Large $v_{1,3}$} ;
\draw [ fill={rgb,255:red,240; green,255; blue,122} ] (6.25,2.75) circle (0.5cm) node {\Large $v_{2,1}$} ;
\draw [ fill={rgb,255:red,255; green,159; blue,220} ] (8.75,2.75) circle (0.5cm) node {\Large $v_{2,2}$} ;
\draw [ fill={rgb,255:red,88; green,249; blue,246} ] (11.25,2.75) circle (0.5cm) node {\Large $v_{2,3}$} ;
\draw [ fill={rgb,255:red,88; green,249; blue,246} ] (13.75,5.25) circle (0.5cm) node {\Large $v_{1,4}$} ;
\draw [ fill={rgb,255:red,240; green,255; blue,122} ] (13.75,2.75) circle (0.5cm) node {\Large $v_{2,4}$} ;
\draw [short] (6.75,5.25) -- (8.25,5.25);
\draw [short] (9.25,5.25) -- (10.75,5.25);
\draw [short] (13.75,4.75) -- (13.75,3.25);
\draw [short] (11.75,5.25) -- (13.25,5.25);
\draw [short] (11.75,2.75) -- (13.25,2.75);
\draw [->, >=Stealth] (12.5,2.75) -- (12.25,2.75) ;
\draw [short] (10.75,2.75) -- (9.25,2.75);
\draw [short] (8.25,2.75) -- (6.75,2.75);
\draw [->, >=Stealth] (6.75,5.25) -- (7.5,5.25);
\draw [->, >=Stealth] (9.25,5.25) -- (10,5.25);
\draw [->, >=Stealth] (12.25,5.25) -- (12.5,5.25);
\draw [->, >=Stealth] (13.75,4.5) -- (13.75,3.75);
\draw [->, >=Stealth] (10.25,2.75) -- (9.75,2.75);
\draw [->, >=Stealth] (7.75,2.75) -- (7.25,2.75);
\node [font=\LARGE] at (7.4,5.75) {$1$};
\node [font=\LARGE] at (9.9,5.75) {$1$};
\node [font=\LARGE] at (12.4,5.75) {$1$};
\node [font=\LARGE] at (14.75,4) {$1$};
\node [font=\LARGE] at (9.6,3.15) {$-1$};
\node [font=\LARGE] at (7.1,3.15) {$-1$};
\node [font=\LARGE] at (12.1,3.15) {$-1$};
\draw [line width=1pt, short] (5.5,4) -- (14.25,4);
\draw [line width=1pt, short] (7.25,4) -- (8.5,7);
\draw [line width=1pt, short] (9.75,4) -- (11,7);
\draw [line width=1pt, short] (12.25,4) -- (13.5,7);
\draw [line width=1pt, dashed] (7.25,4) -- (8.5,1);
\draw [line width=1pt, dashed] (9.75,4) -- (11,1);
\draw [line width=1pt, dashed] (12.25,4) -- (13.5,1);
\end{circuitikz}}}
\subcaptionbox{Weights of edges of $M^*_{2,4}$ under $\tilde{\beta}$ and the $\tilde{h}_{\tilde{\gamma},v_{1,1}}(\tilde{\beta},v_{i,j})$ values (in bold).}[0.3\textwidth]{
\resizebox{0.3\textwidth}{!}{%
\begin{circuitikz}
\draw [ fill={rgb,255:red,88; green,249; blue,246} ] (6.25,5.25) circle (0.5cm) node {\Large $v_{1,1}$} ;
\draw [ fill={rgb,255:red,255; green,159; blue,220} ] (8.75,5.25) circle (0.5cm) node {\Large $v_{1,2}$} ;
\draw [ fill={rgb,255:red,240; green,255; blue,122} ] (11.25,5.25) circle (0.5cm) node {\Large $v_{1,3}$} ;
\draw [ fill={rgb,255:red,255; green,159; blue,220} ] (6.25,2.75) circle (0.5cm) node {\Large $v_{2,1}$} ;
\draw [ fill={rgb,255:red,240; green,255; blue,122} ] (8.75,2.75) circle (0.5cm) node {\Large $v_{2,2}$} ;
\draw [ fill={rgb,255:red,88; green,249; blue,246} ] (11.25,2.75) circle (0.5cm) node {\Large $v_{2,3}$} ;
\draw [ fill={rgb,255:red,88; green,249; blue,246} ] (13.75,5.25) circle (0.5cm) node {\Large $v_{1,4}$} ;
\draw [ fill={rgb,255:red,255; green,159; blue,220} ] (13.75,2.75) circle (0.5cm) node {\Large $v_{2,4}$} ;
\draw [short] (6.75,5.25) -- (8.25,5.25);
\draw [short] (9.25,5.25) -- (10.75,5.25);
\draw [short] (13.75,4.75) -- (13.75,3.25);
\draw [short] (11.75,5.25) -- (13.25,5.25);
\draw [short] (11.75,2.75) -- (13.25,2.75);
\draw [->, >=Stealth] (12.5,2.75) -- (12.25,2.75) ;
\draw [short] (10.75,2.75) -- (9.25,2.75);
\draw [short] (8.25,2.75) -- (6.75,2.75);
\draw [->, >=Stealth] (6.75,5.25) -- (7.5,5.25);
\draw [->, >=Stealth] (9.25,5.25) -- (10,5.25);
\draw [->, >=Stealth] (12.25,5.25) -- (12.5,5.25);
\draw [->, >=Stealth] (13.75,4.5) -- (13.75,3.75);
\draw [->, >=Stealth] (10.25,2.75) -- (9.75,2.75);
\draw [->, >=Stealth] (7.75,2.75) -- (7.25,2.75);
\node [font=\LARGE] at (7.4,5.75) {$-1$};
\node [font=\LARGE] at (9.9,5.75) {$-1$};
\node [font=\LARGE] at (12.4,5.75) {$-1$};
\node [font=\LARGE] at (14.75,4) {$-1$};
\node [font=\LARGE] at (9.8,3.25) {$1$};
\node [font=\LARGE] at (7.3,3.25) {$1$};
\node [font=\LARGE] at (12.2,3.25) {$1$};
\draw [line width=1pt, dashed] (5.5,4) -- (14.25,4);
\draw [line width=1pt, dashed] (7.25,4) -- (8.5,7);
\draw [line width=1pt, dashed] (9.75,4) -- (11,7);
\draw [line width=1pt, dashed] (12.25,4) -- (13.5,7);
\draw [line width=1pt, short] (7.25,4) -- (8.5,1);
\draw [line width=1pt, short] (9.75,4) -- (11,1);
\draw [line width=1pt, short] (12.25,4) -- (13.5,1);
\node [font=\LARGE] at (6.25,6.25) {\textbf{0}};
\node [font=\LARGE] at (8.75,6.25) {$\mathbf{-2}$};
\node [font=\LARGE] at (11.25,6.25) {$\mathbf{-4}$};
\node [font=\LARGE] at (13.75,6.25) {$\mathbf{-6}$};
\node [font=\LARGE] at (13.75,1.85) {$\mathbf{-8}$};
\node [font=\LARGE] at (11.25,1.85) {$\mathbf{-6}$};
\node [font=\LARGE] at (8.75,1.85) {$\mathbf{-4}$};
\node [font=\LARGE] at (6.25,1.85) {$\mathbf{-2}$};
   \end{circuitikz}}}
\subcaptionbox{Graph of the  $ |\tilde{h}_{\tilde{\gamma},v_{1,1}}(\tilde{\beta},v_{i,j})|$ values as $v_{i,j}$ travels along the directed cycle $K^*$, viewed as a Dyck path with maximum area.}[0.3\textwidth]{
\resizebox{0.3\textwidth}{!}{%
    \begin{circuitikz}
     \draw [line width=1pt, short] (0,1) -- (9.5,1);
     \draw [line width=1pt, short] (1,0) -- (1,6);
     \draw [line width=1pt, short] (1,1) -- (5,5) -- (9,1);
     \draw [line width=1pt, short] (2,.7) -- (2,1.3);
     \draw [line width=1pt, short] (3,.7) -- (3,1.3);
     \draw [line width=1pt, short] (4,.7) -- (4,1.3);
     \draw [line width=1pt, short] (5,.7) -- (5,1.3);
     \draw [line width=1pt, short] (6,.7) -- (6,1.3);
     \draw [line width=1pt, short] (7,.7) -- (7,1.3);
     \draw [line width=1pt, short] (8,.7) -- (8,1.3);
     \draw [line width=1pt, short] (9,.7) -- (9,1.3);
     \draw [line width=1pt, short] (.7,2) -- (1.3,2);
     \draw [line width=1pt, short] (.7,3) -- (1.3,3);
     \draw [line width=1pt, short] (.7,4) -- (1.3,4);
     \draw [line width=1pt, short] (.7,5) -- (1.3,5);
     \node [font=\LARGE] at (10,1) {$v_{i,j}$};
     \node [font=\LARGE] at (2,.5) {$v_{1,2}$};
     \node [font=\LARGE] at (3,.5) {$v_{1,3}$};
     \node [font=\LARGE] at (4,.5) {$v_{1,4}$};
     \node [font=\LARGE] at (5,.5) {$v_{2,4}$};
     \node [font=\LARGE] at (6,.5) {$v_{2,3}$};
     \node [font=\LARGE] at (7,.5) {$v_{2,2}$};
     \node [font=\LARGE] at (8,.5) {$v_{2,1}$};
     \node [font=\LARGE] at (9,.5) {$v_{1,1}$};
     \node [font=\LARGE] at (.5,2) {2};
     \node [font=\LARGE] at (.5,3) {4};
     \node [font=\LARGE] at (.5,4) {6};
     \node [font=\LARGE] at (.5,5) {8};
     \node [font=\LARGE] at (1,6.6) {$|\tilde{h}_{\tilde{\gamma},v_{1,1}}(\tilde{\beta},v_{i,j})|$};
     \draw [ fill={rgb,255:red,0; green,0; blue,0} ] (1,1) circle (0.25cm) ; 
     \draw [ fill={rgb,255:red,0; green,0; blue,0} ] (2,2) circle (0.25cm) ; 
     \draw [ fill={rgb,255:red,0; green,0; blue,0} ] (3,3) circle (0.25cm) ; 
     \draw [ fill={rgb,255:red,0; green,0; blue,0} ] (4,4) circle (0.25cm) ; 
     \draw [ fill={rgb,255:red,0; green,0; blue,0} ] (5,5) circle (0.25cm) ; 
     \draw [ fill={rgb,255:red,0; green,0; blue,0} ] (6,4) circle (0.25cm) ; 
     \draw [ fill={rgb,255:red,0; green,0; blue,0} ] (7,3) circle (0.25cm) ; 
     \draw [ fill={rgb,255:red,0; green,0; blue,0} ] (8,2) circle (0.25cm) ; 
     \draw [ fill={rgb,255:red,0; green,0; blue,0} ] (9,1) circle (0.25cm) ; 
     \end{circuitikz}
     }}%
    \caption{3-colorings $\tilde{\gamma}$, $\tilde{\beta}$, and the $\tilde{h}_{\tilde{\gamma},v_{1,1}}(\tilde{\beta},v_{i,j})$ values in $M^*_{2,4}$.}
    \label{fig:M24}
\end{figure}

\end{proof}

While the proof of Theorem~\ref{thm:shortestpath} uses the two MV assignments of $M_{2,n}$ where $\alpha_{1,n}$ and $\alpha_{2,1}$ are the only flippable faces, it also holds, by symmetry, for the other MV assignments that have only two flippable faces (i.e., the four degree-2 vertices in $\OFG(M_{2,n})$). However, in order to obtain the maximum $|h_{\gamma,v_{1,1}}(\beta,v_{i,j})|$ values for the MV assignments whose only flippable faces are $\alpha_{1,1}$ and $\alpha_{2,n}$, one needs to let the base vertex $u$ in Lemma~\ref{lem:3colshortestpath} be $u=v_{2,1}$ instead of $u=v_{1,1}$. Nonetheless, these four degree-2 vertices in $\OFG(M_{2,n})$ are the only MV assignments that, with their opposites, give the maximal-area Dyck paths, as illustrated  in the Figure~\ref{fig:M24} example.

The following corollary results from Lemma \ref{thm:oppositeexist}, Theorem \ref{thm:distbound}, and Theorem \ref{thm:shortestpath}. 
\begin{corollary}
    The diameter of $\OFG(M_{2,n})$ is $\lceil\frac{n^2}{2}\rceil$. 
\end{corollary}

We note that the proof of Theorem~\ref{thm:shortestpath} does not generalize to $M_{m,n}$ for $m>2$. One reason is because in order for the graph of the relative height function $h$ to be like a Dyck path, $M^*_{m,n}$ needs to have a Hamilton cycle. In general, $M^*_{m,n}$ will not have a Hamilton cycle when $m$ and $n$ are both odd, and thus a different argument would need to be made for those cases. However, even for cases like $M^*_{3,4}$, which does have a Hamilton cycle, we have not been able to find a locally valid MV assignment whose corresponding 3-coloring, together with its MV opposite, give a relative height function $h$ whose graph is the maximum-area Dyck path. It seems like the geometry of the Miura-ori crease pattern (and the nuances of 3-coloring grid graphs) limits the maximum value that $h$ can attain, and more insight into this subtlety would be needed to generalize our proof.



\section{Conclusion}\label{sec:conclusion}
We have found several properties of the origami flip graph of the Miura-ori. In $\OFG(M_{2,n})$, we enumerated the number of vertices and edges, proved its diameter, and proved numerous properties about the degree sequence. For the general case, we showed that $\OFG(M_{m,n})$ has exactly four vertices of degree 2 and found bounds on the diameter. This leads us to several new open problems, such as the ones below.
\begin{itemize}
    \item Can we generalize results about $\OFG(M_{2,n})$ to larger cases of $m$? For instance, can we find the diameter of an $m \times n$ Miura-ori, or say anything about the degree sequence?
    \item Can we develop a more straight-forward approach for finding the actual polynomials that describe the number of vertices of $\OFG(M_{2,n})$ with degree $d$?
\end{itemize}

Lastly, we note that as this paper was being revised and prepared for publication, a new paper appeared on the arXiv, \cite{Gupta}, which provides some answers to the above questions. In particular, \cite{Gupta} enumerates the number of degree-4 and degree-5 vertices in $\OFG(M_{m,n})$, gives an alternate proof of our Theorem~\ref{thm:4verticesdeg2}, and finds equations for the diameter of $\OFG(M_{m,n})$ when $m=3$, $m=4$, and $m=n$ that match the bound given in our Theorem~\ref{thm:distbound}.


\section*{Acknowledgements}
This work was supported by NSF grant DMS-2149647, and the second author was supported in part by NSF grant DMS-2347000. The authors thank Mathematical Staircase, Inc. for procuring funding for this research at the 2024 MathILy-EST REU. The authors also thank Jacob Paltrowitz for his proof of Theorem \ref{thm:2xnedges}, as well as the numerous people who provided helpful support and commentary.


\bibliography{refs}
\bibliographystyle{abbrv}

\end{document}